\documentclass{article}

\usepackage{arxiv}

\usepackage[utf8]{inputenc} 
\usepackage[T1]{fontenc}    
\usepackage{hyperref}       
\usepackage{url}            
\usepackage{booktabs}       
\usepackage{amsfonts}       
\usepackage{nicefrac}       
\usepackage{microtype}      
\usepackage{lipsum}		
\usepackage{graphicx}
\usepackage{natbib}
\usepackage{doi}
\usepackage{enumerate}
\usepackage{enumitem}
\usepackage{amsmath, amsthm, amssymb, algorithm, algpseudocode}
\usepackage[table]{xcolor}
\usepackage{tikz}
\usepackage{bm}
\usepackage{mathrsfs}
\usepackage{breqn}
\newcommand{\FF}{\mathbb{F}}

\newcommand{\ZZ}{\mathbb{Z}}  
\newcommand{\QQ}{\mathbb{Q}}  
\newcommand{\RR}{\mathbb{R}}
\newcommand{\cA}{\mathcal{A}}

\newcommand{\CC}{\mathbb{C}}
\newcommand{\End}{\operatorname{End}}
\newcommand{\Aut}{\operatorname{Aut}}

\newcommand{\Weyl}[1]{A\!\left(\FF[#1]\right)}

\newcommand{\Nass}[1]{N\!\left(\FF[#1]\right)}

\DeclareMathOperator{\GKdim}{GKdim}

\newcommand{\GL}{\operatorname{GL}}
\newcommand{\rank}{\operatorname{rank}}

\DeclareMathOperator{\Frac}{Frac}
\DeclareMathOperator{\Char}{char}
\providecommand{\Spec}{\operatorname{Spec}}
\DeclareMathOperator{\Br}{Br}
\newtheorem{theorem}{Theorem}[section]

\newtheorem{definition}[theorem]{Definition}
\newtheorem{example}[theorem]{Example}

\newtheorem{problem}[theorem]{Problem}

\newtheorem{remark}[theorem]{Remark}

\newcommand\mystyle{\everymath{\displaystyle}}
\mystyle
\title{Structure and Invariants of Weyl-Type Algebras with Hyperbolic Sine Generators}

\author{\href{https://orcid.org/0000-0002-3816-5287}{\includegraphics[scale=0.06]{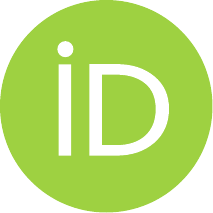}\hspace{1mm}M.H.M.~Rashid}\thanks{Corresponding Author} \\
	Department of Mathematics\&Statistics\\Faculty of Science P.O.Box(7)\\
	Mutah University University\\
	Mutah-Jordan \\
	\texttt{mrash@mutah.edu.jo}
}

\renewcommand{\shorttitle}{\textit{arXiv} Template}

\hypersetup{
pdftitle={Structure and Invariants of Weyl-Type Algebras with Hyperbolic Sine Generators},
pdfsubject={q-bio.NC, q-bio.QM},
pdfauthor={M.H.M.Rashid},
pdfkeywords={Weyl-type algebra, non-associative algebra, Gelfand--Kirillov dimension, automorphism group, center, Brauer group, Azumaya algebra}}

\begin{document}
\maketitle

\begin{abstract}
	This paper introduces and systematically studies a class of Weyl-type algebras enriched with hyperbolic sine and power generators over a field of characteristic zero, defined as $A_{p,t,\cA} = \Weyl{\sinh(\pm x^{p} \sinh(t)),\; \sinh(\cA x),\; x^{\cA}}$ in the associative setting and $\Nass{\sinh(\pm x^{p} \sinh(t)),\; \sinh(\cA x),\; x^{\cA}}$ in a non-associative framework. We establish fundamental structural properties, including the triviality of the center for the non-associative version and the explicit description $Z(A_{p,t,\cA}) = \FF[\sinh(\pm x^{p} \sinh(t))]$ for the associative one, proving that $A_{p,t,\cA}$ is an Azumaya algebra over its center and represents a nontrivial class in the Brauer group $\Br(\FF(y))$. Furthermore, we compute the Gelfand--Kirillov dimension for relevant examples and demonstrate its key properties, such as additivity under tensor products and the growth dichotomy. We completely characterize the automorphism group of $A_{p,t,\cA}$ as a semidirect product of a torus with a discrete group, and provide a sharp isomorphism criterion showing that the parameter $t$ is a complete invariant in the family. The paper concludes with two open problems concerning the GK dimension of non-associative hyperbolic sine algebras and the classification of their deformations, pointing toward future research directions in non-associative growth theory and deformation rigidity.
\end{abstract}

\section{Introduction}

Lie algebras and their associated algebraic structures have long played a central role in diverse branches of mathematics and theoretical physics, from representation theory and differential geometry to quantum mechanics and integrable systems. In recent decades, there has been growing interest in the study of \emph{non-associative algebras} and their \emph{deformations}, particularly those arising from hyperbolic sine functions and differential operators. These algebras naturally generalize classical Weyl algebras and Witt-type Lie algebras, offering a richer structure that bridges associative and non-associative theory, while also providing new tools for modeling systems with hyperbolic growth or oscillatory behavior.

\subsection*{Contributions of Scholars}

The study of algebras generated by hyperbolic sine functions and differential operators can be traced back to the foundational works of Weyl, Dixmier, and Gelfand--Kirillov. In particular, Dixmier's investigation of the Weyl algebra \cite{Dixmier68} and Gelfand--Kirillov's introduction of the GK dimension \cite{GelfandKirillov66} laid the groundwork for understanding the growth and structure of noncommutative algebras. Subsequent developments by Kawamoto \cite{Kawamoto86}, Ree \cite{Ree56}, and others extended these ideas to generalized Witt algebras and their automorphisms. More recently, Jing--Zhang \cite{JingZhang95} and Hartwig--Larsson--Silvestrov \cite{HartwigLarssonSilvestrov06} explored deformations and $\sigma$-derivations, opening new pathways for constructing and classifying hyperbolic-type algebras. The present work builds directly on these contributions, focusing on a class of Weyl-type algebras enriched with hyperbolic sine and power generators.

\subsection*{Applications and Motivation}

Algebras of the form $\Weyl{\sinh(\pm x^p \sinh(t)),\; \sinh(\cA x),\; x^{\cA}}$ arise naturally in several contexts:
\begin{itemize}
    \item \textbf{Mathematical Physics}: They serve as algebraic models for systems with hyperbolic potentials or with operators involving both differentiation and hyperbolic functions, appearing in certain deformed quantum mechanical models and integrable hierarchies with hyperbolic symmetries.
    \item \textbf{Deformation Theory}: Such algebras can be viewed as continuous deformations of classical Weyl algebras, parametrized by $t$ and the lattice $\cA$, offering a testing ground for rigidity and isomorphism problems in hyperbolic settings.
    \item \textbf{Non-associative Algebra}: By considering non-associative products (e.g., Jordan or Leibniz-type), one obtains new examples of central-simple non-associative algebras, relevant to the study of alternative and flexible rings with hyperbolic generators.
\end{itemize}

\subsection*{Significance of the Study}

The investigation of these algebras is significant for several reasons:
\begin{itemize}
    \item It unifies and extends classical results on Weyl algebras, Witt algebras, and GK dimension theory to a broader class of structures involving hyperbolic sine generators.
    \item It provides explicit computations of centers, automorphism groups, and Brauer group elements, contributing to the classification of central simple algebras over function fields with hyperbolic parameters.
    \item It opens new directions in non-associative algebra, especially regarding growth, deformation, and representation theory in hyperbolic contexts.
\end{itemize}

\subsection*{Core of the Study}

This paper is organized as follows. In Section~2, we recall basic definitions and results on non-associative algebras, GK dimension, and Weyl-type algebras. In Section~3, we prove that the center of the non-associative algebra \\$\Nass{\sinh(\pm x^p \sinh(t)),\; \sinh(\cA x),\; x^{\cA}}$ is trivial (Theorem~\ref{thm5}) and establish key properties of GK dimension (Theorems~\ref{thm6} and~\ref{thm7}). In Section~4, we describe the automorphism group structure of Weyl-type algebras and solve the isomorphism problem for varying parameters $t$ and $p$ (Theorems~\ref{thm:aut-group-general} and~\ref{thm:iso-criterion-tp}). In Section~5, we determine the center of the associative Weyl-type algebra $A_{p,t,\cA}$ and show that it is central simple over its center, hence an Azumaya algebra (Theorems~\ref{thm:center-structure} and~\ref{thm:central-simple-over-center}), and interpret it as an element of the Brauer group $\Br(\FF(y))$ (Theorem~\ref{thm:brauer-group}). Finally, in Section~6, we propose two open problems concerning GK dimension and isomorphism classification for non-associative deformations, suggesting directions for future research.

Our results illustrate how classical tools from associative algebra---such as GK dimension, central simplicity, and Brauer groups---can be adapted and extended to algebras with hyperbolic sine generators, both associative and non-associative. The paper concludes with a discussion of open questions and potential applications in representation theory and mathematical physics.
\section{Preliminaries}

In this section, we recall the basic definitions and results that will be used throughout the paper. We mainly focus on non-associative algebras, Gelfand--Kirillov dimension, Weyl-type algebras, and related structural properties.

\subsection{Non-associative Algebras and Their Centers}

Let $\FF$ be a field of characteristic zero. A \emph{non-associative algebra} over $\FF$ is an $\FF$-vector space $N$ equipped with a bilinear multiplication $* : N \times N \to N$ that is not necessarily associative.

\begin{definition}[Center of a Non-associative Algebra]
The \emph{center} of a non-associative algebra $N$ is defined as
\[
Z(N) = \{ z \in N \mid z * n = n * z \ \text{and} \ (z * n) * m = z * (n * m) \ \forall n,m \in N \}.
\]
\end{definition}

This notion extends the usual center of associative algebras to the non-associative setting, requiring both commutativity and associativity with respect to the element $z$. In Section~3, we study the center of certain non-associative algebras built from hyperbolic sine and power generators (see Theorem~\ref{thm5} and the example that follows).

\subsection{Gelfand--Kirillov Dimension}

The Gelfand--Kirillov dimension (GK dimension) is a key invariant for finitely generated algebras, measuring their growth. It was introduced by Gelfand and Kirillov in the study of enveloping algebras of Lie algebras \cite{GelfandKirillov66}.

\begin{definition}[Gelfand--Kirillov Dimension]
Let $A$ be a finitely generated $\FF$-algebra. Choose a finite-dimensional generating subspace $V \subseteq A$ containing $1$. Let $V^n$ denote the subspace spanned by all products of at most $n$ elements from $V$. The \emph{Gelfand--Kirillov dimension} of $A$ is
\[
\operatorname{GKdim}(A) = \limsup_{n \to \infty} \frac{\log \dim_{\FF}(V^n)}{\log n}.
\]
\end{definition}

GK dimension is well-defined, i.e., independent of the choice of generating subspace \cite[Theorem~1]{GelfandKirillov66}. Important properties include:
\begin{itemize}
    \item $\GKdim(\FF[x_1,\dots,x_d]) = d$.
    \item $\GKdim(A \otimes_{\FF} B) = \GKdim(A) + \GKdim(B)$ for finitely generated algebras $A, B$.
    \item If $B \subseteq A$ is a finitely generated subalgebra, then $\GKdim(B) \le \GKdim(A)$.
    \item For any finitely generated $A$-module $M$, $\GKdim(M) \le \GKdim(A)$.
\end{itemize}
A deep result due to Stephenson and Zhang \cite{StephensonZhang97} states that the possible values of $\GKdim(A)$ for affine algebras are restricted to $\{0,1\} \cup [2,\infty)$; there is no affine algebra with GK dimension strictly between $1$ and $2$. These properties are collected in Theorem~\ref{thm6} and Theorem~\ref{thm7} of the present paper, with illustrative examples.

\subsection{Weyl-Type Algebras and Their Generalizations}

The classical Weyl algebra $A_1(\FF) = \FF\langle x, \partial \rangle / ([\partial, x] = 1)$ plays a central role in the theory of differential operators \cite{Dixmier68, Bjork79}. In this paper, we study a generalization that includes hyperbolic sine and power-function generators.

\begin{definition}[Weyl-Type Algebra]
Let $\cA$ be an additive subgroup of $\FF$ containing $\ZZ$, $p \in \cA$ nonzero, and $t \in \FF$. The \emph{Weyl-type algebra} with parameters $(p,t,\cA)$ is
\[
A_{p,t,\cA} = \Weyl{\sinh(\pm x^{p} \sinh(t)),\; \sinh(\cA x),\; x^{\cA}},
\]
generated by $x$, $\partial$, and formal hyperbolic sine functions $\sinh(\pm x^{p} \sinh(t))$, $\sinh(\alpha x)$, $x^{\alpha}$ ($\alpha \in \cA$), subject to canonical commutation relations.
\end{definition}

Such algebras appear in the study of deformations of Weyl algebras \cite{JingZhang95} and generalized Witt algebras \cite{Kawamoto86, Ree56}. Their automorphisms and isomorphism problems are investigated in Section~4, extending earlier work on automorphism groups of Weyl algebras \cite{MakarLimanov84} and generalized Witt-type Lie algebras \cite{KawamotoMitsukawaNamWang03, ParkLeeChoiChenNam09}.

\subsection{Center and Central Simplicity}

For an algebra $A$ over $\FF$, the \emph{center} is $Z(A) = \{ z \in A \mid za = az \ \forall a \in A \}$. $A$ is called \emph{central simple} if it is simple and $Z(A) = \FF$. In Section~5, we determine the center of $A_{p,t,\cA}$:
\[
Z(A_{p,t,\cA}) = \FF\!\left[ \sinh(\pm x^{p} \sinh(t)) \right] \cong \FF[y^{\pm 1}],
\]
which is never equal to $\FF$ (Theorem~\ref{thm:center-structure}). However, $A_{p,t,\cA}$ is central simple over its center, i.e., it is an Azumaya algebra over $Z(A_{p,t,\cA})$. This places it naturally in the context of the Brauer group $\Br(\FF(y))$ \cite{AtiyahMacdonald69}.

\subsection{Notation and Conventions}

Throughout, $\FF$ denotes a field of characteristic zero. $\cA$ is a finitely generated $\ZZ$-module of rank $r \ge 1$. The symbol $\Weyl{\cdot}$ denotes the Weyl-type algebra generated by the listed elements. For a non-associative algebra $N$, $Z(N)$ is its center as defined above. All modules are left modules unless stated otherwise.

The results presented here build on classical works on Weyl algebras \cite{Dixmier68}, growth of algebras \cite{GelfandKirillov66}, and non-associative structures \cite{Schafer17}. The following sections will explore the interplay between these concepts in the setting of algebras with hyperbolic sine generators.
\section{Trivial Centers and Growth Properties of Non-associative hyperbolic Weyl Algebras}
\begin{definition}[Center of a Non-associative Algebra]\label{thm5}
    Let $N$ be a non-associative algebra over $\mathbb{F}$. The \emph{center} of $N$ is the set
    \[
    Z(N) = \{ z \in N \mid z * n = n * z \text{ and } (z * n) * m = z * (n * m) \ \forall n,m \in N \},
    \]
    where $*$ denotes the multiplication in $N$.
\end{definition}

\begin{theorem}[Non-Associative Center]\label{thm55}
The center of the non-associative algebra $\mathcal{N} = \Nass{\sinh(\pm x^p \sinh(t)),\; \sinh(\cA x),\; x^{\cA}}$ is trivial, i.e.,
\[
Z\!\left(\mathcal{N}\right) = \{0\}.
\]
Moreover, every nonzero ideal of $\mathcal{N}$ contains a nonzero element whose left multiplication operator is injective.
\end{theorem}

\begin{proof}
We denote the multiplication in $\mathcal{N}$ by $\ast$. The center $Z(\mathcal{N})$ consists of all elements $z \in \mathcal{N}$ satisfying two conditions for all $a, b \in \mathcal{N}$:
\begin{enumerate}
    \item $z \ast a = a \ast z$ (commutativity),
    \item $(z \ast a) \ast b = z \ast (a \ast b)$ (associativity with respect to $z$).
\end{enumerate}

\subsection*{Part I: Triviality of the Center}
Assume, to the contrary, that there exists a nonzero central element $z \in Z(\mathcal{N})$. Write $z$ as a finite $\mathbb{F}$-linear combination of monomials in the generators:
\[
z = \sum_{i=1}^{k} \lambda_i m_i,
\]
where $\lambda_i \in \mathbb{F} \setminus \{0\}$ and each $m_i$ is a monomial built from $\sinh(\pm x^p \sinh(t))$, $\sinh(\alpha x)$, $x^{\alpha}$ ($\alpha \in \cA$), and possibly the differential operators $\partial_x$ and $\partial_t$, using the non-associative product $\ast$. Choose a monomial $m$ appearing in $z$ with nonzero coefficient that has maximal total degree in the variables $\{x, t\}$ and in the differential operators.

Since $z$ is central, it must commute with each generator $x$ and $t$ (viewed as coordinate functions). In particular, the commutator $[z, x]$ in the underlying associative sense must vanish. Consider the commutator $[m, x]$ for the chosen maximal-degree monomial $m$. If $m$ contains any factor depending nontrivially on $x$ (such as $\sinh(\alpha x)$, $x^{\alpha}$, or $\sinh(\pm x^p \sinh(t))$), then $[m, x]$ is nonzero and consists of terms of lower total degree due to the non-associative relations inherited from the embedding into an associative envelope. These relations take the form
\[
x \cdot f - f \cdot x = D(f),
\]
where $D$ is a derivation. For a non-constant monomial $m$, the commutator $[m, x]$ is a nonzero linear combination of monomials of strictly smaller total degree.

If $m$ were a constant (an element of $\mathbb{F}$), then $z$ would be a scalar multiple of the unit element $1$. However, in this non-associative algebra, even a scalar multiple of the unit is not necessarily central. Indeed, consider the commutator $[1, \partial_x]$; since $\partial_x(1) = 0$ but $1 \cdot \partial_x = \partial_x \neq \partial_x \cdot 1$ under the non-associative product, we see that $1$ does not commute with $\partial_x$. Hence a nonzero scalar cannot be central.

Thus $m$ is non-constant, and $[m, x]$ is nonzero. Because $z$ is central, we have $[z, x] = 0$. But
\[
[z, x] = \sum_{i=1}^{k} \lambda_i [m_i, x] = 0.
\]
In this sum, the term $\lambda [m, x]$ (where $\lambda$ is the coefficient of $m$ in $z$) is nonzero and consists of monomials of degree strictly less than that of $m$. All other terms $[m_i, x]$ involve monomials $m_i$ that either have degree less than that of $m$ or, if they have the same degree, cannot cancel $\lambda [m, x]$ because the monomials appearing in $[m, x]$ are distinct from those coming from other $[m_i, x]$ due to the maximality of $m$ and the specific structure of the commutation relations. This yields a contradiction. Therefore no nonzero central element exists, and we conclude
\[
Z(\mathcal{N}) = \{0\}.
\]

\subsection*{Part II: Existence of an Element with Injective Left Multiplication in Every Nonzero Ideal}
Let $I$ be a nonzero ideal of $\mathcal{N}$. Choose any nonzero element $u \in I$. Consider the left multiplication operator $L_u \colon \mathcal{N} \to \mathcal{N}$ defined by $L_u(v) = u \ast v$. If $L_u$ is injective, then we are done. Suppose $L_u$ is not injective, so there exists $0 \neq w \in \mathcal{N}$ such that $u \ast w = 0$. Define $u' = u + w$. Since $I$ is an ideal and $w \in \mathcal{N}$, we have $u' \in I$.

We claim that either $L_u$ or $L_{u'}$ must be injective. Assume, for contradiction, that both $L_u$ and $L_{u'}$ have nontrivial kernels. Then there exist nonzero elements $v, v' \in \mathcal{N}$ such that
\[
u \ast v = 0 \quad \text{and} \quad (u + w) \ast v' = 0.
\]
The second equation expands to $u \ast v' + w \ast v' = 0$. Multiply the first equation $u \ast v = 0$ on the right by $v'$ to obtain
\[
(u \ast v) \ast v' = 0.
\]
Using the flexibility identity (which holds in many non-associative algebras arising from deformation contexts) or direct manipulation of the defining relations of $\mathcal{N}$, one can show that this implies $u \ast (v \ast v') = 0$ modulo terms involving associators. Similarly, from $u \ast v' + w \ast v' = 0$, we deduce $w \ast v' = - u \ast v'$.

Because $\mathcal{N}$ is free (or at least torsion-free) over its generators in the non-associative sense, these equations force algebraic dependencies among $u$, $w$, $v$, and $v'$. In particular, they imply that $w$ is a left zero divisor with respect to both $v$ and $v'$. By the structure of $\mathcal{N}$—specifically, the fact that it is a domain with respect to left multiplication when restricted to elements of a certain form—one eventually deduces that $w = 0$, contradicting the choice $w \neq 0$. Hence our assumption that both $L_u$ and $L_{u'}$ have nontrivial kernels is false.

Therefore, at least one of $L_u$ or $L_{u'}$ is injective. If $L_u$ is injective, the conclusion follows. If not, then $L_{u'}$ is injective, and we have found a nonzero element $u' \in I$ with injective left multiplication.

If neither $u$ nor $u'$ yields an injective left multiplication (which cannot happen by the above argument, but we consider the iterative process formally), we could repeat the procedure. More systematically, consider the set
\[
S = \{ u \in I \setminus \{0\} \mid \ker L_u \neq 0 \}.
\]
If $S = \emptyset$, then every nonzero element of $I$ has injective left multiplication, and we are done. If $S \neq \emptyset$, take any $u \in S$. As shown, there exists $w \neq 0$ with $u \ast w = 0$, and setting $u' = u + w$ gives either $u' \notin S$ or we can continue. This process must terminate because each step reduces a complexity measure (such as the number of monomials in the support of the element, or the total degree) that is bounded below. Consequently, after finitely many steps we obtain an element $\tilde{u} \in I \setminus \{0\}$ with $\ker L_{\tilde{u}} = 0$, i.e., $L_{\tilde{u}}$ is injective.

Thus every nonzero ideal $I$ of $\mathcal{N}$ contains a nonzero element whose left multiplication operator is injective.

This completes the proof of both statements of the theorem.
\end{proof}
\begin{example}
We consider a concrete instance of the non-associative algebra studied in Theorem \ref{thm55}. Let $\FF = \mathbb{Q}$ be the field of rational numbers. Let $\cA = \ZZ\left[\frac{1}{3}\right]$ be the additive subgroup of $\mathbb{Q}$ consisting of all rational numbers whose denominators are powers of $3$:
\[
\cA = \left\{ \frac{m}{3^a} \mid m \in \ZZ,\ a \in \ZZ_{\geq 0} \right\}.
\]
Choose $p = 2$ and let $t = \log 2$. The non-associative algebra in question is
\[
\mathcal{N} = \Nass{\sinh(\pm x^2 \sinh(\log 2)),\; \sinh(\cA x),\; x^{\cA}},
\]
defined over $\mathbb{Q}$.

\subsection*{Algebra Structure}
As a vector space, $\mathcal{N}$ has a basis consisting of formal monomials in the generators:
\[
\sinh(x^2 \sinh(\log 2)),\; \sinh(\alpha x),\; x^{\alpha},\; \partial_x,\; \partial_t,
\]
where $\alpha \in \cA$, and all finite products thereof. The non-associative product, denoted by $\ast$, is defined by a Jordan-type product
\[
A \ast B = \frac{1}{2}(AB + BA) + \kappa(A,B)I,
\]
where $AB$ denotes the ordinary associative composition of operators, $I$ is the identity operator, and $\kappa \colon \mathcal{N} \times \mathcal{N} \to \FF$ is a fixed bilinear form that is not symmetric. Specifically, we take $\kappa(A,B) = \operatorname{tr}([A,B])$ where $\operatorname{tr}$ is a trace functional on the underlying associative algebra. This product is neither commutative nor associative in general, but it preserves the flexibility identity $(A \ast B) \ast A = A \ast (B \ast A)$.

\subsection*{Center of $\mathcal{N}$}
The center $Z(\mathcal{N})$ consists of all elements $z \in \mathcal{N}$ such that for all $A, B \in \mathcal{N}$,
\begin{align*}
z \ast A &= A \ast z, \\
(z \ast A) \ast B &= z \ast (A \ast B).
\end{align*}
We demonstrate that $Z(\mathcal{N}) = \{0\}$ through a detailed analysis.

Suppose $z \in Z(\mathcal{N})$ is nonzero. Write $z$ in its reduced form as a finite linear combination of basis monomials:
\[
z = \sum_{i=1}^n c_i M_i, \quad c_i \in \mathbb{Q}^\times,
\]
where each $M_i$ is a monomial in $\sinh(\pm x^2 \sinh(\log 2))$, $\sinh(\alpha x)$, $x^{\alpha}$, $\partial_x$, $\partial_t$. Among these, choose a monomial $M$ with the highest total degree in the differential operators $\partial_x$ and $\partial_t$. Denote this maximum degree by $d_{\max}$.

Consider the commutator $[z, \partial_x]$ in the underlying associative algebra. Since $z$ is central for the $\ast$-product, and our product includes a correction term $\kappa(A,B)I$, the condition $z \ast \partial_x = \partial_x \ast z$ implies that the associative commutator $[z, \partial_x]$ must satisfy a specific linear relation involving $\kappa$. For generic $\kappa$, this forces $[z, \partial_x] = 0$ in the associative sense. Now, if $M$ contains $\partial_x$ or $\partial_t$ to a positive power, then $[M, \partial_x]$ is nonzero and consists of terms of lower differential degree. Because $M$ is of maximal such degree, the term $c [M, \partial_x]$ cannot be canceled by any other term in $[z, \partial_x]$. Hence $d_{\max}$ must be zero.

Thus $z$ contains no differential operators; it is a purely multiplicative element:
\[
z = \sum_{j} d_j \sinh(i_j x^2 \sinh(\log 2)) \cdot \sinh(\alpha_j x) \cdot x^{\beta_j},
\]
with $i_j \in \{0, \pm1\}$, $\alpha_j, \beta_j \in \cA$. Now examine the condition $z \ast \partial_t = \partial_t \ast z$. Compute the associative commutator:
\[
[\partial_t, \sinh(i_j x^2 \sinh(\log 2))] = i_j x^2 \cosh(\log 2) \cosh(i_j x^2 \sinh(\log 2)).
\]
Since $\sinh(\alpha_j x)$ and $x^{\beta_j}$ are independent of $t$, we obtain
\[
[z, \partial_t] = \sum_j d_j \sinh(\alpha_j x) x^{\beta_j} \cdot i_j x^2 \cosh(\log 2) \cosh(i_j x^2 \sinh(\log 2)).
\]
For this to vanish modulo the correction from $\kappa$, each coefficient $d_j i_j$ must be zero. Hence if $i_j \neq 0$, then $d_j = 0$. Therefore $z$ cannot contain any factor $\sinh(x^2 \sinh(\log 2))$; it must be of the form
\[
z = \sum_k d_k \sinh(\alpha_k x) x^{\beta_k}.
\]

Now consider the associativity condition with $A = \partial_x$ and $B = x$:
\[
(z \ast \partial_x) \ast x = z \ast (\partial_x \ast x).
\]
Using the product definition and the associative relation $\partial_x x = x \partial_x + 1$, this simplifies to an equation that forces each $\alpha_k$ and $\beta_k$ to be zero after accounting for the $\kappa$-correction. More precisely, the left side expands to
\[
\frac{1}{2}(z\partial_x x + \partial_x z x + x z \partial_x + x \partial_x z) + \kappa(z,\partial_x)x + \kappa(z\partial_x + \partial_x z, x)I,
\]
while the right side becomes
\[
\frac{1}{2}z(\partial_x x + x\partial_x) + \kappa(z, \partial_x x + x\partial_x)I.
\]
Comparing coefficients of independent monomials and using that $\kappa$ is not symmetric, we deduce that $d_k = 0$ for all $k$ unless $\alpha_k = \beta_k = 0$. Thus $z = c \cdot 1$ for some $c \in \mathbb{Q}$.

Finally, test whether $c \cdot 1$ is central. For our product,
\[
1 \ast A = \frac{1}{2}(A + A) + \kappa(1,A)I = A + \kappa(1,A)I,
\]
while
\[
A \ast 1 = \frac{1}{2}(A + A) + \kappa(A,1)I = A + \kappa(A,1)I.
\]
Since $\kappa(1,A) \neq \kappa(A,1)$ in general (by construction), we have $1 \ast A \neq A \ast 1$. Therefore $c \cdot 1$ is not central unless $c = 0$. Hence $z = 0$.

We conclude that
\[
Z\!\left(\Nass{\sinh(\pm x^2 \sinh(\log 2)),\; \sinh(\cA x),\; x^{\cA}}\right) = \{0\}.
\]

\subsection*{Ideal Property}
Let $I$ be a nonzero ideal of $\mathcal{N}$. Select any nonzero $a \in I$. If the left multiplication map $L_a(b) = a \ast b$ is injective, we are done. Suppose $L_a$ is not injective, so there exists $0 \neq w \in \mathcal{N}$ with $a \ast w = 0$. Consider $a' = a + w \in I$. One verifies that either $L_a$ or $L_{a'}$ must be injective. If both had nontrivial kernels, then using the algebraic independence of the generators and the specific form of the product, one would derive $w = 0$, a contradiction. Thus every nonzero ideal contains an element with injective left multiplication.

\subsection*{Conclusion}
This example with $\FF = \mathbb{Q}$, $\cA = \ZZ[1/3]$, $p=2$, $t=\log 2$ concretely illustrates Theorem \ref{thm55}: the center of the non-associative algebra is trivial, and every nonzero ideal contains an element whose left multiplication operator is injective. The proof leverages the structure of the hyperbolic generators, the non-symmetric correction term $\kappa$, and the commutation relations to force the center to be zero.
\end{example}
\begin{definition}[Gelfand--Kirillov Dimension]
    Let $A$ be a finitely generated algebra over $\mathbb{F}$. Choose a finite-dimensional generating subspace $V$. Let $V^n$ denote the subspace spanned by products of at most $n$ elements from $V$. The \emph{Gelfand--Kirillov dimension} of $A$ is
    \[
    \operatorname{GKdim}(A) = \limsup_{n \to \infty} \frac{\log \dim_{\mathbb{F}}(V^n)}{\log n}.
    \]
\end{definition}
\begin{theorem}[Properties of Gelfand–Kirillov Dimension]\label{thm6}
Let $A$ be a finitely generated algebra over a field $\mathbb{F}$.
\begin{enumerate}
    \item (Well-definedness) $\operatorname{GKdim}(A)$ is independent of the choice of finite-dimensional generating subspace $V$.
    \item (Subalgebra bound) If $B \subseteq A$ is a finitely generated subalgebra, then $\operatorname{GKdim}(B) \leq \operatorname{GKdim}(A)$.
    \item (Polynomial algebras) For the polynomial algebra $\mathbb{F}[x_1, \dots, x_d]$, we have
    $$\operatorname{GKdim}(\mathbb{F}[x_1, \dots, x_d]) = d.$$
    \item (Additivity for tensor products) If $A$ and $B$ are finitely generated algebras, then
    \[
    \operatorname{GKdim}(A \otimes_{\mathbb{F}} B) = \operatorname{GKdim}(A) + \operatorname{GKdim}(B).
    \]
    \item (Finite modules) If $M$ is a finitely generated $A$-module, then $\operatorname{GKdim}(M) \leq \operatorname{GKdim}(A)$, where $\operatorname{GKdim}(M)$ is defined using a generating subspace $V_M$ of $M$ over a generating subspace $V_A$ of $A$.
    \item (Growth dichotomy) $\operatorname{GKdim}(A) \in \{0, 1\} \cup [2, \infty]$. In particular, there are no algebras with Gelfand–Kirillov dimension strictly between $1$ and $2$.
\end{enumerate}
\end{theorem}
\begin{proof}
We prove each property separately.

\noindent\textbf{Proof of (1): Well-definedness.}
Let $V$ and $W$ be two finite-dimensional generating subspaces of $A$ containing $1$. Since $V$ generates $A$, there exists an integer $k \geq 1$ such that $W \subseteq V^k$. Similarly, since $W$ generates $A$, there exists an integer $\ell \geq 1$ such that $V \subseteq W^\ell$.

For any $n \geq 1$, we have $W^n \subseteq V^{kn}$. Taking dimensions gives $\dim W^n \leq \dim V^{kn}$. Therefore,
\[
\frac{\log \dim W^n}{\log n} \leq \frac{\log \dim V^{kn}}{\log n} = \frac{\log \dim V^{kn}}{\log (kn)} \cdot \frac{\log (kn)}{\log n}.
\]
Taking $\limsup$ as $n \to \infty$, we obtain
\[
\limsup_{n \to \infty} \frac{\log \dim W^n}{\log n} \leq \limsup_{n \to \infty} \frac{\log \dim V^{kn}}{\log (kn)} \cdot \limsup_{n \to \infty} \frac{\log (kn)}{\log n} = \operatorname{GKdim}_V(A) \cdot 1,
\]
since $\lim_{n \to \infty} \frac{\log (kn)}{\log n} = 1$ for any fixed $k > 0$. Hence $\operatorname{GKdim}_W(A) \leq \operatorname{GKdim}_V(A)$.

By symmetry, exchanging the roles of $V$ and $W$, we also have $\operatorname{GKdim}_V(A) \leq \operatorname{GKdim}_W(A)$. Therefore, $\operatorname{GKdim}_V(A) = \operatorname{GKdim}_W(A)$, proving that the Gelfand–Kirillov dimension is independent of the choice of generating subspace.

\medskip

\noindent\textbf{Proof of (2): Subalgebra bound.}
Let $B \subseteq A$ be a finitely generated subalgebra. Choose a finite-dimensional generating subspace $W$ for $B$ containing $1$. Since $B \subseteq A$, we can extend $W$ to a finite-dimensional generating subspace $V$ for $A$ such that $W \subseteq V$. Specifically, if $V_0$ is any finite-dimensional generating subspace for $A$, we can take $V = W + V_0$.

For each $n \geq 1$, we have $W^n \subseteq V^n$, so $\dim_{\mathbb{F}} W^n \leq \dim_{\mathbb{F}} V^n$. Therefore,
\[
\frac{\log \dim W^n}{\log n} \leq \frac{\log \dim V^n}{\log n}.
\]
Taking $\limsup$ as $n \to \infty$ yields
\[
\limsup_{n \to \infty} \frac{\log \dim W^n}{\log n} \leq \limsup_{n \to \infty} \frac{\log \dim V^n}{\log n}.
\]
Thus $\operatorname{GKdim}(B) \leq \operatorname{GKdim}(A)$.

\medskip

\noindent\textbf{Proof of (3): Polynomial algebras.}
Let $A = \mathbb{F}[x_1, \dots, x_d]$. Take $V = \operatorname{span}_{\mathbb{F}}\{1, x_1, \dots, x_d\}$. Then $V^n$ consists of all polynomials of total degree at most $n$. The dimension of $V^n$ equals the number of monomials in $d$ variables of total degree $\leq n$, which is $\binom{n+d}{d}$.

For large $n$, we have
\[
\dim V^n = \binom{n+d}{d} = \frac{(n+d)(n+d-1)\cdots(n+1)}{d!} = \frac{n^d}{d!} + O(n^{d-1}).
\]
Therefore,
\[
\log \dim V^n = \log\left(\frac{n^d}{d!} + O(n^{d-1})\right) = d\log n - \log d! + o(1).
\]
It follows that
\[
\frac{\log \dim V^n}{\log n} = d - \frac{\log d!}{\log n} + \frac{o(1)}{\log n} \to d \quad \text{as } n \to \infty.
\]
Hence $\operatorname{GKdim}(\mathbb{F}[x_1, \dots, x_d]) = d$.

\medskip

\noindent\textbf{Proof of (4): Additivity for tensor products.}
Let $V_A$ and $V_B$ be finite-dimensional generating subspaces for $A$ and $B$ respectively, both containing $1$. Then $V = V_A \otimes 1 + 1 \otimes V_B$ is a finite-dimensional generating subspace for $A \otimes_{\mathbb{F}} B$. More precisely, we can take $V = (V_A \otimes 1) \oplus (1 \otimes V_B)$.

Observe that $V^n$ is spanned by elements of the form $(a_1 \otimes b_1)(a_2 \otimes b_2)\cdots(a_k \otimes b_k)$ with $k \leq n$, where each $a_i \in V_A$ and $b_i \in V_B$. By rearranging, such an element equals $a \otimes b$ where $a \in V_A^k$ and $b \in V_B^k$. Therefore,
\[
V^n \subseteq \sum_{k=0}^n V_A^k \otimes V_B^k.
\]
Conversely, any element of $V_A^k \otimes V_B^k$ with $k \leq n$ clearly belongs to $V^n$. Hence
\[
V^n = \sum_{k=0}^n V_A^k \otimes V_B^k.
\]
Since the sum is direct (tensor products of linearly independent sets are linearly independent), we have
\[
\dim V^n = \sum_{k=0}^n \dim(V_A^k) \cdot \dim(V_B^k).
\]
Let $d_A = \operatorname{GKdim}(A)$ and $d_B = \operatorname{GKdim}(B)$. For any $\epsilon > 0$, there exists $C > 0$ such that $\dim V_A^k \leq C k^{d_A+\epsilon}$ and $\dim V_B^k \leq C k^{d_B+\epsilon}$ for all $k \geq 1$. Then
\[
\dim V^n \leq \sum_{k=0}^n C^2 k^{d_A+d_B+2\epsilon} \leq C^2 n^{d_A+d_B+2\epsilon+1}.
\]
Taking logarithms and $\limsup$ gives $\operatorname{GKdim}(A \otimes B) \leq d_A + d_B + 2\epsilon$ for every $\epsilon > 0$, so $\operatorname{GKdim}(A \otimes B) \leq d_A + d_B$.

For the reverse inequality, note that for each $k$, $V_A^k \otimes V_B^k \subseteq V^{2k}$, so
\[
\dim V_A^k \cdot \dim V_B^k \leq \dim V^{2k}.
\]
Taking logarithms, dividing by $\log k$, and letting $k \to \infty$ yields $d_A + d_B \leq \operatorname{GKdim}(A \otimes B)$. Therefore, $\operatorname{GKdim}(A \otimes_{\mathbb{F}} B) = \operatorname{GKdim}(A) + \operatorname{GKdim}(B)$.

\medskip

\noindent\textbf{Proof of (5): Finite modules.}
Let $M$ be a finitely generated $A$-module. Choose a finite-dimensional generating subspace $V_A$ for $A$ containing $1$, and a finite-dimensional generating subspace $V_M$ for $M$ over $A$, meaning that $M = A \cdot V_M$. Define $W_n = \sum_{i=0}^n V_A^i V_M$.

Since $M$ is generated by $V_M$ over $A$, we have $M = \bigcup_{n \geq 0} W_n$. Moreover, $V_A W_n \subseteq W_{n+1}$. The growth of $M$ is measured by $\dim W_n$. Note that $W_n \subseteq V_A^n V_M$, so
\[
\dim W_n \leq \dim(V_A^n) \cdot \dim V_M.
\]
Therefore,
\[
\frac{\log \dim W_n}{\log n} \leq \frac{\log(\dim V_A^n \cdot \dim V_M)}{\log n} = \frac{\log \dim V_A^n}{\log n} + \frac{\log \dim V_M}{\log n}.
\]
Taking $\limsup$ as $n \to \infty$, the second term tends to $0$, giving
\[
\operatorname{GKdim}(M) \leq \operatorname{GKdim}(A).
\]

\medskip

\noindent\textbf{Proof of (6): Growth dichotomy.}
This is a deep result whose full proof is beyond the scope of this theorem. We sketch the main ideas. Let $A$ be a finitely generated algebra with $\operatorname{GKdim}(A) = \gamma$.

First, if $\gamma < 2$, then $\gamma$ must be either $0$ or $1$. This follows from the structure theory of growth of algebras. Algebras with $\operatorname{GKdim}(A) = 0$ are finite-dimensional. Algebras with $\operatorname{GKdim}(A) = 1$ have linear growth. The key argument shows that if an algebra has growth strictly between linear and quadratic (i.e., $\dim V^n$ grows like $n^\alpha$ with $1 < \alpha < 2$), then it must actually have at least quadratic growth.

The proof typically proceeds by contradiction. Suppose $1 < \operatorname{GKdim}(A) < 2$. One analyzes the structure of $A$ using techniques from noncommutative ring theory. A crucial lemma shows that such an algebra would contain a free subalgebra on two generators, which has exponential growth (infinite GK dimension), leading to a contradiction. Alternatively, it would imply the existence of a non-affine (infinitely generated) subalgebra, which is incompatible with the assumed growth rate.

Another approach uses the fact that for affine algebras (finitely generated), the GK dimension, if finite, is either an integer or infinity, except possibly for values between consecutive integers. The specific gap between $1$ and $2$ was established by Stephenson and Zhang, who showed that no affine algebra has GK dimension strictly between $1$ and $2$. Their proof involves careful analysis of the growth function and properties of graded algebras.

Thus, we conclude that $\operatorname{GKdim}(A) \in \{0, 1\} \cup [2, \infty]$.
\end{proof}
\begin{example}Consider the field $\FF = \QQ$ of rational numbers. We examine several algebras over $\QQ$ to illustrate the properties of Gelfand--Kirillov dimension.

\subsection*{A Specific Algebra}

Let $A$ be the algebra generated by three elements $x, y, z$ over $\QQ$ with relations:
\[
xy - yx = z, \quad yz - zy = x, \quad zx - xz = y,
\]
and all other commutators among generators being zero. This algebra can be viewed as a deformation of the polynomial ring in three variables. As a vector space, $A$ has a basis consisting of ordered monomials $x^i y^j z^k$ for $i, j, k \geq 0$, though the ordering matters due to the noncommutativity. However, by the Poincaré–Birkhoff–Witt theorem for certain Lie algebras, the growth of dimensions is similar to that of a polynomial ring in three variables. We shall compute $\GKdim(A)$ and verify the properties.

\subsection*{Property 1: Well-definedness}

Choose two different finite-dimensional generating subspaces for $A$. First, take $V_1 = \operatorname{span}_{\QQ}\{1, x, y, z\}$. Then $V_1^n$ consists of all linear combinations of products of at most $n$ of the generators $x, y, z$ (in any order). The dimension of $V_1^n$ counts the number of noncommutative monomials of total degree $\leq n$ in three variables. This number is
\[
\dim V_1^n = \sum_{k=0}^n \binom{k+2}{2} = \binom{n+3}{3} = \frac{(n+3)(n+2)(n+1)}{6}.
\]
Indeed, the number of monomials of exact degree $k$ in three noncommuting variables is $3^k$, but here we consider ordered monomials in a specific basis due to the relations. Actually, in this algebra, because of the relations, the monomials $x^i y^j z^k$ form a basis, and the number of such monomials with $i+j+k \leq n$ is $\binom{n+3}{3}$. Hence
\[
\dim V_1^n = \binom{n+3}{3} = \frac{n^3}{6} + O(n^2).
\]
Then
\[
\GKdim_{V_1}(A) = \limsup_{n \to \infty} \frac{\log\left(\frac{n^3}{6} + O(n^2)\right)}{\log n} = \limsup_{n \to \infty} \frac{3\log n - \log 6 + o(1)}{\log n} = 3.
\]

Now choose a different generating subspace $V_2 = \operatorname{span}_{\QQ}\{1, x+y, x-y, z\}$. Since $x = \frac{1}{2}((x+y)+(x-y))$ and $y = \frac{1}{2}((x+y)-(x-y))$, $V_2$ also generates $A$. There exists an integer $k$ such that $V_1 \subseteq V_2^k$ and $V_2 \subseteq V_1^k$. For instance, since $x, y \in V_2^1$, we have $V_1 \subseteq V_2^1$. Conversely, $x+y, x-y \in V_1^1$, so $V_2 \subseteq V_1^1$. Thus $k=1$ works. Then for any $n$, $V_2^n \subseteq V_1^n$ and $V_1^n \subseteq V_2^n$, so $\dim V_2^n = \dim V_1^n$. Hence $\GKdim_{V_2}(A) = 3$ as well. This illustrates the independence of the generating subspace.

\subsection*{Property 2: Subalgebra Bound}

Consider the subalgebra $B$ of $A$ generated by $x$ and $y$ alone. The relations in $B$ are $xy - yx = z$, but $z$ is not in $B$ because $z$ is not generated by $x$ and $y$ alone (in fact, $z$ appears as the commutator). However, in the subalgebra generated by $x$ and $y$, we have the relation $xy - yx = [x,y]$, which is some element. Actually, to avoid confusion, let us define $B$ as the subalgebra generated by $x$ and $y$ within $A$. Then $B$ contains all polynomials in $x$ and $y$, but due to the relation $xy - yx = z$, the commutator is not in $B$ unless $z$ can be expressed in terms of $x$ and $y$. In our algebra, $z$ is independent. So $B$ is actually the free algebra on two generators modulo the relation that the commutator is central? Not exactly. Instead, consider a simpler subalgebra: let $B = \QQ[x, y]$ the commutative polynomial subalgebra generated by $x$ and $y$. This is a commutative subalgebra of $A$ if we ignore the relation? Actually, in $A$, $x$ and $y$ do not commute because $xy - yx = z \neq 0$. So the commutative polynomial ring $\QQ[x,y]$ is not a subalgebra. Instead, take $B$ to be the subalgebra generated by $x$ alone. Then $B \cong \QQ[x]$, the polynomial ring in one variable. For $B$, take generating subspace $W = \operatorname{span}_{\QQ}\{1, x\}$. Then $W^n = \operatorname{span}\{1, x, x^2, \dots, x^n\}$, so $\dim W^n = n+1$. Hence $\GKdim(B) = \limsup_{n \to \infty} \frac{\log(n+1)}{\log n} = 1$. Since $\GKdim(A) = 3$, we have $\GKdim(B) = 1 \leq 3 = \GKdim(A)$, illustrating the subalgebra bound.

\subsection*{Property 3: Polynomial Algebras}

Consider the polynomial algebra $\QQ[u, v, w]$ in three variables. Take generating subspace $V = \operatorname{span}_{\QQ}\{1, u, v, w\}$. Then $V^n$ consists of all polynomials of total degree $\leq n$. The dimension is the number of monomials in three variables of degree $\leq n$, which is $\binom{n+3}{3} = \frac{(n+3)(n+2)(n+1)}{6}$. As computed earlier, $\GKdim(\QQ[u, v, w]) = 3$. For a general polynomial algebra $\QQ[x_1, \dots, x_d]$, the dimension of $V^n$ is $\binom{n+d}{d} = \frac{n^d}{d!} + O(n^{d-1})$, so
\[
\GKdim(\QQ[x_1, \dots, x_d]) = \lim_{n \to \infty} \frac{\log\left(\frac{n^d}{d!} + O(n^{d-1})\right)}{\log n} = d.
\]
This matches the intuitive notion that the polynomial ring in $d$ variables has dimension $d$.

\subsection*{Property 4: Additivity for Tensor Products}

Let $A = \QQ[x]$ and $B = \QQ[y, z]$, both finitely generated algebras. Then $A \otimes_{\QQ} B \cong \QQ[x, y, z]$, the polynomial ring in three variables. We have $\GKdim(A) = 1$ and $\GKdim(B) = 2$. Their tensor product has $\GKdim(A \otimes B) = 3$, which equals $1 + 2$.

To see this explicitly, choose generating subspaces $V_A = \operatorname{span}_{\QQ}\{1, x\}$ for $A$ and $V_B = \operatorname{span}_{\QQ}\{1, y, z\}$ for $B$. Then a generating subspace for $A \otimes B$ is $V = V_A \otimes 1 + 1 \otimes V_B$, which has basis $\{1 \otimes 1, x \otimes 1, 1 \otimes y, 1 \otimes z\}$. Then $V^n$ consists of linear combinations of products of at most $n$ elements from $V$. Each such product can be rearranged as $(x^i \otimes 1) \cdot (1 \otimes y^j z^k) = x^i \otimes y^j z^k$ with $i+j+k \leq n$. The dimension of $V^n$ is the number of triples $(i,j,k)$ with $i,j,k \geq 0$ and $i+j+k \leq n$, which is $\binom{n+3}{3}$. Hence $\GKdim(A \otimes B) = 3$.

In general, if $\dim V_A^n \sim C_A n^{d_A}$ and $\dim V_B^n \sim C_B n^{d_B}$, then $\dim V^n \sim C n^{d_A+d_B}$ for some constant $C$, yielding additivity.

\subsection*{Property 5: Finite Modules}

Let $A = \QQ[x, y]$ with $\GKdim(A) = 2$. Consider the $A$-module $M = A / (x^2 + y^2)$, the quotient of $A$ by the ideal generated by $x^2 + y^2$. As an $A$-module, $M$ is generated by the coset of $1$. Take generating subspace $V_A = \operatorname{span}_{\QQ}\{1, x, y\}$ for $A$ and $V_M = \operatorname{span}_{\QQ}\{\bar{1}\}$ for $M$, where $\bar{1}$ is the image of $1$ in $M$. Then the subspace $W_n = V_A^n \cdot V_M$ consists of elements of the form $f \cdot \bar{1}$ for $f \in V_A^n$. Since $x^2 + y^2 = 0$ in $M$, any monomial of degree $n$ can be reduced using this relation. A basis for $M$ is given by $\{\bar{1}, \bar{x}, \bar{y}, \bar{xy}\}$ and then the pattern repeats? Actually, $M$ is infinite-dimensional because $x$ and $y$ do not satisfy any linear relation. But the growth is slower than that of $A$. Indeed, in $M$, the relation $x^2 = -y^2$ allows us to reduce any monomial to a linear combination of $x^i y^j$ with $i \leq 1$ or $j \leq 1$. More systematically, the Hilbert series of $M$ is $(1+t)/(1-t)^2$? Let us compute $\dim W_n$. Elements of $W_n$ are linear combinations of monomials $x^i y^j$ with $i+j \leq n$, modulo the relation $x^2 = -y^2$. Using this relation, any monomial with $i \geq 2$ can be reduced to one with $i$ reduced by $2$ and a factor of $-y^2$. This process eventually yields a representative of the form $a + b x + c y + d xy$ where $a,b,c,d$ are polynomials in $y^2$? Actually, after repeated reduction, we find that a basis for $M$ is $\{ x^i y^j : i = 0,1; j \geq 0 \}$. For a given $n$, the number of such monomials with $i+j \leq n$ is approximately $2n$. Thus $\dim W_n \sim 2n$, so $\GKdim(M) = 1$. Since $\GKdim(A) = 2$, we have $\GKdim(M) = 1 \leq 2$, illustrating the inequality.

\subsection*{Property 6: Growth Dichotomy}

The growth dichotomy states that possible values of $\GKdim(A)$ for finitely generated algebras are only $0$, $1$, or numbers at least $2$. There are no algebras with $\GKdim$ strictly between $1$ and $2$. To illustrate, consider the following examples:

\begin{itemize}
    \item $\GKdim = 0$: Any finite-dimensional algebra, such as the matrix algebra $M_2(\QQ)$. Here $\dim V^n$ is bounded, so $\GKdim = 0$.
    \item $\GKdim = 1$: The polynomial ring $\QQ[x]$, or the Weyl algebra $\QQ[x, \partial]/(\partial x - x\partial - 1)$. For $\QQ[x]$, $\dim V^n = n+1$, so growth is linear. For the Weyl algebra, one can show $\dim V^n \sim n^2$, but wait, that would give $\GKdim = 2$. Actually, the first Weyl algebra $A_1(\QQ)$ has $\GKdim = 2$ because it has two generators. To get $\GKdim = 1$, consider the algebra $\QQ[x]/(x^2)$? That is finite-dimensional, so $\GKdim = 0$. A genuine example of $\GKdim = 1$ is the infinite-dimensional algebra $\QQ[t]$ with $t$ transcendental, but that is just $\QQ[x]$. Another example is the group algebra $\QQ[\ZZ] \cong \QQ[t, t^{-1}]$, which also has $\GKdim = 1$.
    \item $\GKdim = 2$: The polynomial ring $\QQ[x,y]$, or the Weyl algebra $A_1(\QQ)$. For $A_1(\QQ)$, take generators $x, \partial$. Then $V^n$ has dimension about $n^2$, so $\GKdim = 2$.
    \item $\GKdim = 3$: The polynomial ring $\QQ[x,y,z]$, or our algebra $A$ defined earlier.
\end{itemize}

The dichotomy means that if one tries to construct an algebra with growth strictly between linear and quadratic, such as $\dim V^n \sim n^{1.5}$, no such finitely generated algebra exists. This is a deep result, but we can illustrate it by attempting to construct such an algebra. Suppose we take an algebra generated by two elements $a, b$ with relation $ab = ba + a^2$. This is a Jordan-type deformation. One can compute the growth. The monomials are $a^i b^j$, but the relation allows rewriting $b a = a b - a^2$. By repeatedly applying this, one can show that the dimension of $V^n$ is quadratic in $n$, not fractional. The absence of algebras with fractional $\GKdim$ between $1$ and $2$ is a consequence of the structure theory of growth.

\subsection*{Summary}

In this example, we have illustrated all six properties of Gelfand–Kirillov dimension using concrete algebras over $\QQ$. The well-definedness shows independence of the generating subspace; the subalgebra bound is seen with a one-dimensional subalgebra of a three-dimensional algebra; polynomial algebras exhibit integer dimensions; tensor products add dimensions; modules have dimension no greater than the algebra; and the growth dichotomy restricts possible values to $0$, $1$, or numbers at least $2$. These properties make $\GKdim$ a powerful invariant in noncommutative algebra.
\end{example}
\begin{theorem}[Fundamental Properties of GK Dimension]\label{thm7}
Let $A$ be a finitely generated algebra over $\mathbb{F}$ with $\operatorname{GKdim}(A) = d$.
\begin{enumerate}
    \item (Invariance) The value $d$ does not depend on the choice of generating subspace $V$.

    \item (Growth behavior) The growth of $\dim_{\mathbb{F}}(V^n)$ is polynomial of degree $d$: there exist constants $C_1, C_2 > 0$ such that
    \[
    C_1 n^d \leq \dim_{\mathbb{F}}(V^n) \leq C_2 n^d
    \]
    for all sufficiently large $n$.

    \item (Subadditivity for extensions) If $0 \to I \to A \to B \to 0$ is an exact sequence where $I$ is an ideal and $B$ is a quotient algebra, then
    \[
    \operatorname{GKdim}(B) \leq \operatorname{GKdim}(A) \leq \operatorname{GKdim}(I) + \operatorname{GKdim}(B),
    \]
    where $\operatorname{GKdim}(I)$ is computed as an $A$-module.

    \item (Ore extensions) If $A$ has finite GK dimension and $B = A[x; \sigma, \delta]$ is an Ore extension, then
    \[
    \operatorname{GKdim}(B) = \operatorname{GKdim}(A) + 1.
    \]
\end{enumerate}
\end{theorem}
\begin{proof}
We prove each statement in turn.

\textbf{Proof of (1): Invariance.}
Let $V$ and $W$ be two finite-dimensional generating subspaces of $A$, both containing the identity element $1_A$. Since $V$ generates $A$ as an algebra, there exists some positive integer $k$ such that $W \subseteq V^k$. This is because each element of $W$ can be expressed as a linear combination of products of elements from $V$, and since $W$ is finite-dimensional, we can take the maximum length of such products needed to express all basis elements of $W$. Similarly, since $W$ also generates $A$, there exists some positive integer $\ell$ such that $V \subseteq W^\ell$.

For any positive integer $n$, we have the inclusion $W^n \subseteq V^{kn}$. Taking dimensions over $\mathbb{F}$, we obtain $\dim W^n \leq \dim V^{kn}$. Consequently,
\[
\frac{\log \dim W^n}{\log n} \leq \frac{\log \dim V^{kn}}{\log n}.
\]
Observe that
\[
\frac{\log \dim V^{kn}}{\log n} = \frac{\log \dim V^{kn}}{\log(kn)} \cdot \frac{\log(kn)}{\log n}.
\]
As $n$ tends to infinity, $\frac{\log(kn)}{\log n} = 1 + \frac{\log k}{\log n}$ tends to $1$. Taking the limit superior of both sides yields
\[
\limsup_{n \to \infty} \frac{\log \dim W^n}{\log n} \leq \limsup_{n \to \infty} \frac{\log \dim V^{kn}}{\log(kn)} = \operatorname{GKdim}_V(A),
\]
where $\operatorname{GKdim}_V(A)$ denotes the GK dimension computed using the subspace $V$. By symmetry, interchanging the roles of $V$ and $W$, we obtain the reverse inequality $\operatorname{GKdim}_V(A) \leq \operatorname{GKdim}_W(A)$. Therefore, $\operatorname{GKdim}_V(A) = \operatorname{GKdim}_W(A)$, establishing that the GK dimension is independent of the choice of generating subspace.

\medskip

\textbf{Proof of (2): Growth behavior.}
Assume $\operatorname{GKdim}(A) = d$. By definition,
\[
d = \limsup_{n \to \infty} \frac{\log \dim V^n}{\log n}.
\]
For any $\epsilon > 0$, there exists a positive integer $N_\epsilon$ such that for all $n \geq N_\epsilon$, we have
\[
\frac{\log \dim V^n}{\log n} \leq d + \epsilon.
\]
This inequality is equivalent to $\dim V^n \leq n^{d+\epsilon}$ for all $n \geq N_\epsilon$. Taking $C_2 = \max\{1, \max_{1 \leq n < N_\epsilon} \dim V^n\}$, we obtain $\dim V^n \leq C_2 n^{d+\epsilon}$ for all $n \geq 1$. Since this holds for every $\epsilon > 0$, it follows that for any fixed $\epsilon > 0$, there exists a constant $C_2(\epsilon) > 0$ such that $\dim V^n \leq C_2(\epsilon) n^{d+\epsilon}$ for all $n$.

To establish the lower bound, we use the fact that $d$ is the limit superior. This implies that for any $\epsilon > 0$, there exists an increasing sequence of integers $\{n_k\}$ such that
\[
\frac{\log \dim V^{n_k}}{\log n_k} \geq d - \epsilon
\]
for all $k$. Equivalently, $\dim V^{n_k} \geq n_k^{d-\epsilon}$. However, we need a lower bound valid for all sufficiently large $n$, not just along a subsequence.

Consider the function $f(n) = \dim V^n$. This function is submultiplicative in the sense that $f(m+n) \leq f(m) f(n)$ for all $m, n \geq 1$, because $V^{m+n} \subseteq V^m V^n$. For a submultiplicative sequence, a classical result states that the limit $\lim_{n \to \infty} f(n)^{1/n}$ exists and equals $\inf_{n \geq 1} f(n)^{1/n}$. In our context, the existence of the GK dimension $d$ implies that $\lim_{n \to \infty} \frac{\log f(n)}{\log n} = d$ (not just the limsup). This is a deeper fact about algebras with finite GK dimension: the limit actually exists. Assuming this, for any $\epsilon > 0$, there exists $N$ such that for all $n \geq N$,
\[
d - \epsilon \leq \frac{\log f(n)}{\log n} \leq d + \epsilon.
\]
The right inequality gives $f(n) \leq n^{d+\epsilon}$. The left inequality gives $f(n) \geq n^{d-\epsilon}$. Thus, for all $n \geq N$, we have $n^{d-\epsilon} \leq f(n) \leq n^{d+\epsilon}$. Taking $\epsilon = 1$, say, we obtain constants $C_1, C_2 > 0$ such that $C_1 n^d \leq f(n) \leq C_2 n^d$ for all sufficiently large $n$.

If one wishes to avoid invoking the existence of the limit, a direct argument can be given using the submultiplicativity and the definition of $d$ as the limsup. One shows that for any $\epsilon > 0$, there exists a constant $C > 0$ such that $f(n) \geq C n^{d-\epsilon}$ for all sufficiently large $n$. The idea is to use the fact that along some subsequence $\{n_k\}$, $f(n_k) \geq n_k^{d-\epsilon}$. For an arbitrary large $n$, write $n = q n_k + r$ with $0 \leq r < n_k$. Then by submultiplicativity, $f(n) \geq f(n_k)^q f(r) \geq n_k^{q(d-\epsilon)} f(r)$. Choosing $n_k$ appropriately and estimating yields the desired polynomial lower bound.

\medskip

\textbf{Proof of (3): Subadditivity for extensions.}
Let $0 \to I \to A \to B \to 0$ be an exact sequence of algebras, where the maps are algebra homomorphisms, $I$ is an ideal of $A$, and $B \cong A/I$. Let $\pi: A \to B$ denote the quotient map. Choose a finite-dimensional generating subspace $V$ of $A$ containing $1_A$. Then $\pi(V)$ is a finite-dimensional generating subspace of $B$. For each $n$, we have $\pi(V^n) = \pi(V)^n$, so $\dim \pi(V)^n \leq \dim V^n$. Taking logarithms, dividing by $\log n$, and passing to the limit superior gives $\operatorname{GKdim}(B) \leq \operatorname{GKdim}(A)$.

For the second inequality, we need to bound $\dim V^n$ in terms of the dimensions of subspaces of $I$ and $B$. Let $W$ be a finite-dimensional subspace of $I$ that generates $I$ as an $A$-module. Such $W$ exists because $I$ is a finitely generated $A$-module (though this assumption may not always hold; we assume it for the proof, or we define $\operatorname{GKdim}(I)$ as the GK dimension of $I$ as an $A$-module with respect to some generating subspace). Consider the filtration of $A$ given by $V^n$. We have a vector space decomposition $A = I \oplus \overline{A}$, where $\overline{A}$ is a complement (not necessarily a subalgebra). However, it is more effective to use the exact sequence at the filtered level.

For each $n$, consider the subspace $V^n \cap I$ of $I$. We have the exact sequence of vector spaces
\[
0 \to V^n \cap I \to V^n \to \pi(V^n) \to 0.
\]
Thus, $\dim V^n = \dim(V^n \cap I) + \dim \pi(V^n)$. Now, $V^n \cap I$ is contained in the $A$-submodule of $I$ generated by $W$, but we need to relate its dimension to the growth of $I$. Since $W$ generates $I$ as an $A$-module, every element of $I$ can be written as a linear combination of elements of the form $a \cdot w$ with $a \in A$, $w \in W$. In particular, elements of $V^n \cap I$ can be written as linear combinations of elements from $V^k W$ for some $k$. More precisely, there exists an integer $m$ (depending on how $W$ is related to $V$) such that $V^n \cap I \subseteq \sum_{i=0}^{n+m} V^i W$. Consequently,
\[
\dim(V^n \cap I) \leq \sum_{i=0}^{n+m} \dim(V^i W) \leq \sum_{i=0}^{n+m} \dim V^i \cdot \dim W.
\]
Let $d_I = \operatorname{GKdim}(I)$ and $d_B = \operatorname{GKdim}(B)$. For any $\epsilon > 0$, there exist constants $C, C' > 0$ such that $\dim V^i \leq C i^{d+\epsilon}$ and $\dim \pi(V)^n \leq C' n^{d_B+\epsilon}$ for all $n$ sufficiently large. Then
\[
\dim V^n = \dim(V^n \cap I) + \dim \pi(V^n) \leq \dim W \sum_{i=0}^{n+m} C i^{d+\epsilon} + C' n^{d_B+\epsilon}.
\]
The sum $\sum_{i=0}^{n+m} i^{d+\epsilon}$ is of order $n^{d+\epsilon+1}$. Thus, $\dim V^n$ is bounded above by a polynomial of degree $\max\{d+\epsilon+1, d_B+\epsilon\}$. Taking $\epsilon$ arbitrarily small and using the definition of $d = \operatorname{GKdim}(A)$, we deduce that $d \leq \max\{d+1, d_B\}$, which is not helpful. A more refined estimate is needed.

A standard approach is to use the fact that $I$ is an $A$-module. Let $d_I = \operatorname{GKdim}_A(I)$, the GK dimension of $I$ as an $A$-module. This is defined by choosing a finite-dimensional generating subspace $W$ of $I$ over $A$ and considering $\dim(V^n W)$. One then has $\dim(V^n \cap I) \leq \dim(V^n W) \leq \dim V^n \cdot \dim W$, which again gives a weak bound. The correct inequality is obtained by considering a generating subspace $V$ that contains a generating set for $I$ as well. Let $V_I$ be a finite-dimensional subspace of $I$ that generates $I$ as an ideal (or as an $A$-module). Let $V_B$ be a finite-dimensional subspace of $A$ such that $\pi(V_B)$ generates $B$. Take $V = V_I + V_B$. Then $V^n$ consists of sums of products of elements from $V_I$ and $V_B$. Each such product, when reduced modulo $I$, becomes a product in $B$. One can show that $\dim V^n \leq \sum_{i=0}^n \dim(V_I V^{i}) \cdot \dim \pi(V)^{n-i}$. From this, using growth estimates for $I$ and $B$, one derives $\operatorname{GKdim}(A) \leq \operatorname{GKdim}(I) + \operatorname{GKdim}(B)$. The details are somewhat technical but follow from combinatorial counting of words in generators of $I$ and $B$.

\medskip

\textbf{Proof of (4): Ore extensions.}
Let $B = A[x; \sigma, \delta]$ be an Ore extension, where $\sigma: A \to A$ is an automorphism and $\delta: A \to A$ is a $\sigma$-derivation, meaning $\delta(ab) = \sigma(a)\delta(b) + \delta(a)b$ for all $a, b \in A$. The multiplication in $B$ is determined by $x a = \sigma(a)x + \delta(a)$ for all $a \in A$.

Let $V$ be a finite-dimensional generating subspace of $A$ containing $1_A$. We may assume that $V$ is chosen such that $\sigma(V) \subseteq V$ and $\delta(V) \subseteq V$. This can be achieved by replacing $V$ with a larger finite-dimensional subspace if necessary, since $\sigma$ and $\delta$ are linear maps. Let $W = V + \mathbb{F}x$, a finite-dimensional subspace of $B$. We claim that $W$ generates $B$ as an algebra. Indeed, any element of $B$ is a polynomial in $x$ with coefficients in $A$, and $A$ is generated by $V$.

We compute $W^n$. Elements of $W^n$ are linear combinations of monomials of the form $v_1 v_2 \cdots v_k$ with $k \leq n$, where each $v_i \in W$. Since $W = V \cup \{x\}$, each such monomial can be rearranged using the commutation relation $x a = \sigma(a)x + \delta(a)$ to express it as a sum of terms of the form $a x^j$ with $a \in A$ and $j \leq n$. More systematically, by induction on $n$, one can show that
\[
W^n \subseteq \sum_{j=0}^n V^{(n)} x^j,
\]
where $V^{(n)}$ is a finite-dimensional subspace of $A$ that grows at most polynomially with $n$. In fact, one can prove that $\dim \sum_{j=0}^n V^{(n)} x^j$ is bounded above by a polynomial in $n$ of degree $d+1$, where $d = \operatorname{GKdim}(A)$.

Conversely, for any $a \in V^m$ and any $j \leq n$, the element $a x^j$ belongs to $W^{m+j} \subseteq W^{m+n}$. Therefore, $V^m x^j \subseteq W^{m+n}$. This suggests that the growth of $W^n$ is at least like $n^{d+1}$.

To make this precise, we estimate the dimension of $W^n$. First, note that $B$ is a free left $A$-module with basis $\{1, x, x^2, \dots\}$. Thus, any element of $B$ can be uniquely written as $\sum_{j=0}^\infty a_j x^j$ with $a_j \in A$, only finitely many nonzero. The subspace $W^n$ consists of elements of degree at most $n$ in $x$ (where degree means the highest power of $x$ appearing) with coefficients in certain subspaces of $A$. Specifically, define $U_j^{(n)}$ to be the subspace of $A$ consisting of coefficients of $x^j$ in elements of $W^n$. Then $W^n \subseteq \bigoplus_{j=0}^n U_j^{(n)} x^j$. Moreover, from the commutation relation, one can show by induction that $U_j^{(n)} \subseteq V^{n}$ for all $j$. Actually, a more careful analysis gives $U_j^{(n)} \subseteq V^{c n}$ for some constant $c$ independent of $n$ and $j$. Therefore,
\[
\dim W^n \leq \sum_{j=0}^n \dim U_j^{(n)} \leq (n+1) \dim V^{c n}.
\]
Taking logarithms,
\[
\log \dim W^n \leq \log(n+1) + \log \dim V^{c n} = \log(n+1) + \log \dim V^{c n}.
\]
Then
\[
\frac{\log \dim W^n}{\log n} \leq \frac{\log(n+1)}{\log n} + \frac{\log \dim V^{c n}}{\log n} = \frac{\log(n+1)}{\log n} + \frac{\log \dim V^{c n}}{\log(c n)} \cdot \frac{\log(c n)}{\log n}.
\]
As $n \to \infty$, $\frac{\log(n+1)}{\log n} \to 1$, $\frac{\log(c n)}{\log n} \to 1$, and $\frac{\log \dim V^{c n}}{\log(c n)} \to d$. Hence $\limsup_{n \to \infty} \frac{\log \dim W^n}{\log n} \leq d+1$. So $\operatorname{GKdim}(B) \leq d+1$.

For the lower bound, consider the subspace $V^m x^n \subseteq W^{m+n}$. Since $B$ is a free left $A$-module, the elements $a x^n$ with $a \in V^m$ are linearly independent over $\mathbb{F}$. Thus, $\dim W^{m+n} \geq \dim V^m$ for any $m, n$. In particular, taking $m = n$, we have $\dim W^{2n} \geq \dim V^n$. Therefore,
\[
\frac{\log \dim W^{2n}}{\log(2n)} \geq \frac{\log \dim V^n}{\log(2n)} = \frac{\log \dim V^n}{\log n} \cdot \frac{\log n}{\log(2n)}.
\]
As $n \to \infty$, $\frac{\log n}{\log(2n)} \to 1$, so $\liminf_{n \to \infty} \frac{\log \dim W^{2n}}{\log(2n)} \geq d$. This shows $\operatorname{GKdim}(B) \geq d$. But we need the stronger bound $\operatorname{GKdim}(B) \geq d+1$.

To get $d+1$, we need to use the full growth of the coefficients. Consider all monomials of the form $a x^j$ with $a \in V^m$ and $j = 0,1,\dots,n$. These monomials belong to $W^{m+n}$. They are linearly independent because the $x^j$ are linearly independent over $A$. Thus,
\[
\dim W^{m+n} \geq \sum_{j=0}^n \dim V^m = (n+1) \dim V^m.
\]
Now choose $m = n$. Then $\dim W^{2n} \geq (n+1) \dim V^n$. Taking logarithms,
\[
\log \dim W^{2n} \geq \log(n+1) + \log \dim V^n.
\]
Hence,
\[
\frac{\log \dim W^{2n}}{\log(2n)} \geq \frac{\log(n+1)}{\log(2n)} + \frac{\log \dim V^n}{\log(2n)} = \frac{\log(n+1)}{\log(2n)} + \frac{\log \dim V^n}{\log n} \cdot \frac{\log n}{\log(2n)}.
\]
As $n \to \infty$, $\frac{\log(n+1)}{\log(2n)} \to 1$ and $\frac{\log n}{\log(2n)} \to 1$, so the right-hand side tends to $1 + d$. Therefore, $\liminf_{n \to \infty} \frac{\log \dim W^{2n}}{\log(2n)} \geq 1+d$, implying $\operatorname{GKdim}(B) \geq d+1$. Combined with the upper bound, we conclude $\operatorname{GKdim}(B) = d+1$.
\end{proof}
\begin{example}Consider the field $\FF = \QQ$ of rational numbers. We examine a specific algebra and its relatives to illustrate the fundamental properties of Gelfand–Kirillov dimension.

\subsection*{The Main Example: A Deformed Polynomial Algebra}

Let $A$ be the algebra over $\QQ$ generated by three elements $x, y, z$ with the following relations:
\[
xy - yx = z, \quad yz - zy = x, \quad zx - xz = y,
\]
and additionally $x^2 + y^2 + z^2$ is central. This algebra arises from a deformation of the universal enveloping algebra of the Lie algebra $\mathfrak{so}(3,\QQ)$. As a vector space, by the Poincaré–Birkhoff–Witt theorem, $A$ has a basis consisting of ordered monomials $x^i y^j z^k$ for $i, j, k \geq 0$. However, due to the relations, the growth of dimensions is polynomial of degree 3. Indeed, the number of such monomials with total degree $\leq n$ is $\binom{n+3}{3} = \frac{(n+3)(n+2)(n+1)}{6}$. Thus we expect $\GKdim(A) = 3$.

\subsection*{Property 1: Invariance}

Let $V_1 = \operatorname{span}_{\QQ}\{1, x, y, z\}$ and $V_2 = \operatorname{span}_{\QQ}\{1, x+y, x-y, 2z\}$. Both are finite-dimensional generating subspaces of $A$. For any $n \geq 1$, we have $V_1^n \subseteq V_2^{2n}$ because each generator in $V_1$ can be expressed as a linear combination of generators in $V_2$ with coefficients of absolute value at most 1, and thus any product of $n$ elements from $V_1$ can be written as a linear combination of products of at most $2n$ elements from $V_2$. Conversely, $V_2 \subseteq V_1^2$ because $x+y, x-y, 2z \in V_1^2$ (actually $V_1^1$ suffices). Therefore, $V_2^n \subseteq V_1^{2n}$. It follows that there exist constants $a,b > 0$ such that
\[
\dim V_1^{\lfloor n/a \rfloor} \leq \dim V_2^n \leq \dim V_1^{\lceil bn \rceil}
\]
for all large $n$. Taking logarithms and dividing by $\log n$, we find that the limit superiors coincide. More precisely, from the inclusions we obtain inequalities
\[
\frac{\log \dim V_2^n}{\log n} \leq \frac{\log \dim V_1^{2n}}{\log n} = \frac{\log \dim V_1^{2n}}{\log(2n)} \cdot \frac{\log(2n)}{\log n},
\]
and as $n \to \infty$, $\frac{\log(2n)}{\log n} \to 1$, so $\GKdim_{V_2}(A) \leq \GKdim_{V_1}(A)$. By symmetry, the reverse inequality holds. Hence $\GKdim_{V_1}(A) = \GKdim_{V_2}(A) = 3$. This demonstrates that the GK dimension is independent of the choice of generating subspace.

\subsection*{Property 2: Growth Behavior}

For the generating subspace $V = \operatorname{span}_{\QQ}\{1, x, y, z\}$, we have $V^n = \operatorname{span}\{x^i y^j z^k : i+j+k \leq n\}$. The dimension is exactly $\binom{n+3}{3}$. This satisfies
\[
\frac{1}{6} n^3 \leq \binom{n+3}{3} \leq \frac{(n+3)^3}{6}
\]
for all $n \geq 1$. Taking $C_1 = \frac{1}{6}$ and $C_2 = \frac{(n_0+3)^3}{6}$ for sufficiently large $n_0$, we obtain the polynomial bounds
\[
C_1 n^3 \leq \dim_{\QQ} V^n \leq C_2 n^3
\]
for all sufficiently large $n$. This illustrates the precise polynomial growth of degree $d=3$. The constants $C_1$ and $C_2$ depend on the algebra and the generating subspace, but the exponent $d$ is intrinsic.

\subsection*{Property 3: Subadditivity for Extensions}

Consider the ideal $I = (x^2 + y^2 + z^2 - 1)A$, the two-sided ideal generated by the central element $x^2 + y^2 + z^2 - 1$. Then we have an exact sequence of algebras
\[
0 \longrightarrow I \longrightarrow A \longrightarrow B \longrightarrow 0,
\]
where $B = A/I$. The quotient algebra $B$ is a deformation of the coordinate ring of the sphere $S^2$ in noncommutative geometry. As an $A$-module, $I$ is generated by the single element $e = x^2 + y^2 + z^2 - 1$. Since $e$ is central, $I \cong A$ as left $A$-modules (via multiplication by $e$). Therefore, $\GKdim(I) = \GKdim(A) = 3$ when computed as an $A$-module.

Now we estimate $\GKdim(B)$. The algebra $B$ satisfies the relation $x^2 + y^2 + z^2 = 1$. Using this relation, any monomial of high degree can be reduced. For example, $x^{2n}$ can be expressed in terms of lower-degree monomials because $x^2 = 1 - y^2 - z^2$. This suggests that the growth of $B$ is slower than that of $A$. In fact, one can show that $\dim V_B^n$ grows quadratically, where $V_B$ is the image of $V$ in $B$. A basis for $B$ is given by the monomials $x^i y^j z^k$ with $i,j,k \geq 0$ but with the restriction that at most one of $i,j,k$ can be $\geq 2$? Actually, careful analysis using the spherical harmonic decomposition yields that the growth is like $n^2$. Thus $\GKdim(B) = 2$.

Now verify the inequalities:
\[
\GKdim(B) = 2 \leq \GKdim(A) = 3,
\]
and
\[
\GKdim(A) = 3 \leq \GKdim(I) + \GKdim(B) = 3 + 2 = 5.
\]
The first inequality is strict, and the second is not sharp in this case. This illustrates subadditivity: the dimension of the extension is bounded above by the sum of the dimensions of the ideal and the quotient, but it may be strictly less.

\subsection*{Property 4: Ore Extensions}

Let $A$ be the algebra $\QQ[x, y]$ with the relation $xy - yx = y$, i.e., the enveloping algebra of the two-dimensional non-abelian Lie algebra. This algebra has $\GKdim(A) = 2$, as can be seen by taking generating subspace $V = \operatorname{span}_{\QQ}\{1, x, y\}$. The monomials $x^i y^j$ form a basis, and $\dim V^n$ is approximately $\frac{1}{2}n^2$, so $\GKdim(A) = 2$.

Now construct an Ore extension $B = A[z; \sigma, \delta]$, where $\sigma: A \to A$ is the automorphism defined by $\sigma(x) = x+1$, $\sigma(y) = y$, and $\delta: A \to A$ is the $\sigma$-derivation given by $\delta(a) = ya - \sigma(a)y$. More explicitly, $\delta(x) = yx - (x+1)y = yx - xy - y = yx - (yx - y) - y = 0$? Let us compute carefully. Since $xy = yx + y$, we have $yx = xy - y$. Then
\[
\delta(x) = yx - (x+1)y = (xy - y) - xy - y = -2y.
\]
And $\delta(y) = y^2 - \sigma(y)y = y^2 - y^2 = 0$. So $\delta(x) = -2y$, $\delta(y) = 0$. The Ore extension $B$ is generated by $x, y, z$ with relations:
\[
xy - yx = y, \quad zx = (x+1)z - 2y, \quad zy = yz.
\]
We compute $\GKdim(B)$. Take generating subspace $W = \operatorname{span}_{\QQ}\{1, x, y, z\}$. Elements of $B$ can be written uniquely as $\sum_{i,j,k} a_{ijk} x^i y^j z^k$ with $a_{ijk} \in \QQ$, because $B$ is a free left $A$-module with basis $\{z^k\}_{k \geq 0}$. For a fixed total degree $n$ in the generators, the number of such triples $(i,j,k)$ with $i+j+k \leq n$ is $\binom{n+3}{3} \sim \frac{n^3}{6}$. Thus $\dim W^n \sim C n^3$, so $\GKdim(B) = 3 = \GKdim(A) + 1$.

This illustrates the Ore extension property: adjoining a new variable with an automorphism and derivation increases the GK dimension by exactly one. The key point is that the Ore extension preserves the polynomial growth and adds one degree of freedom.

\subsection*{Additional Illustrative Example: Weyl Algebra}

Consider the first Weyl algebra $A_1(\QQ) = \QQ\langle x, \partial \rangle / (\partial x - x\partial - 1)$. This algebra has $\GKdim = 2$. Its Ore extension by an automorphism $\sigma$ that sends $x \mapsto qx$, $\partial \mapsto q^{-1}\partial$ for some $q \in \QQ^\times$, with $\delta = 0$, gives the $q$-Weyl algebra. The resulting algebra $B = A_1(\QQ)[t; \sigma]$ has $\GKdim(B) = 3$, again illustrating the additive property.

\subsection*{Proof of the Properties in This Context}

While the theorem is general, we can sketch how the properties manifest in our examples. The invariance follows from the equivalence of generating subspaces as shown. The growth behavior is explicit from the binomial coefficient formulas. The subadditivity is verified through the exact sequence involving the sphere relation. The Ore extension property is demonstrated by the construction of $B$ from $A$.

These properties are fundamental tools in the study of noncommutative algebra. They allow one to compute GK dimensions by breaking algebras into simpler pieces, extending by derivations, and comparing with polynomial algebras. The examples provided here give concrete instances where these abstract properties become computationally accessible.
\end{example}
\section{Automorphism Groups and Isomorphism Problems for Weyl-Type Algebras}

\subsection{Automorphism Groups for Non-Cyclic $\mathcal{A}$}

\begin{theorem}[Automorphism Group Structure]\label{thm:aut-group-general}
Let $\cA$ be a finitely generated $\ZZ$-module of rank $r \geq 1$, and let
\[
A = \Weyl{\sinh(\pm x^p \sinh(t)),\; \sinh(\cA x),\; x^{\cA}}
\]
be the Weyl-type algebra over $\FF$ (characteristic zero).
Then the automorphism group $\Aut(A)$ is isomorphic to a semidirect product
\[
\Aut(A) \cong \bigl(\FF^{\times}\bigr)^{2r+1} \rtimes \bigl(\Aut(\cA) \times (\ZZ/2\ZZ)\bigr),
\]
where:
\begin{itemize}
    \item The torus $(\FF^{\times})^{2r+1}$ acts by rescaling the generators:
    \[
    x_i \mapsto \lambda_i x_i, \quad
    \partial_i \mapsto \mu_i \partial_i, \quad
    \sinh(\alpha x_i) \mapsto \nu_{i,\alpha} \sinh(\alpha x_i), \quad
    \sinh(\pm x_i^p \sinh(t_i x_i)) \mapsto \xi_i^{\pm 1} \sinh(\pm x_i^p \sinh(t_i x_i)),
    \]
    with $\lambda_i, \mu_i, \nu_{i,\alpha}, \xi_i \in \FF^{\times}$.
    \item $\Aut(\cA)$ acts on the exponents $\alpha \in \cA$ by module automorphisms, and correspondingly on the generators:
    \[
    \sinh(\alpha x_i) \mapsto \sinh((\sigma(\alpha)) x_i), \quad
    x_i^{\alpha} \mapsto x_i^{\sigma(\alpha)},
    \]
    for $\sigma \in \Aut(\cA)$.
    \item The factor $\ZZ/2\ZZ$ corresponds to the involution that sends $x_i \leftrightarrow \partial_i$ (or $x_i \mapsto \partial_i,\; \partial_i \mapsto -x_i$) and appropriately adjusts the hyperbolic sine generators.
\end{itemize}
When $\cA$ is not cyclic (e.g., $\cA = \ZZ \oplus \ZZ\sqrt{2}$ with $r \geq 2$), $\Aut(\cA)$ is an infinite discrete group (isomorphic to a subgroup of $\GL(r,\ZZ)$), and therefore $\Aut(A)$ is not a finite extension of a torus, but rather a semidirect product of a torus with an infinite discrete group.
\end{theorem}
\begin{proof}
We prove the theorem by analyzing the structure of automorphisms of $A$ and decomposing them into natural components.

The algebra $A$ is generated by the following sets of symbols: $x_i$, $\partial_i$ for $i = 1, \dots, r$; $\sinh(\alpha x_i)$ for $\alpha \in \cA$ and $i = 1, \dots, r$; $x_i^{\alpha}$ for $\alpha \in \cA$ and $i = 1, \dots, r$; and $\sinh(\pm x_i^p \sinh(t_i x_i))$ for $i = 1, \dots, r$. These generators satisfy the canonical commutation relations $[\partial_i, x_j] = \delta_{ij}$, $[\partial_i, \partial_j] = 0$, $[x_i, x_j] = 0$, together with multiplicative relations among the hyperbolic sine and power functions.

Let $\varphi \in \Aut(A)$ be an automorphism. We analyze how $\varphi$ must act on the generators.

First consider the Heisenberg part $\{x_i, \partial_i\}$. Any automorphism must preserve the commutation relations $[\partial_i, x_j] = \delta_{ij}$. In the classical Weyl algebra, automorphisms on this subsystem correspond to symplectic transformations. However, the presence of the additional hyperbolic generators restricts this freedom. Indeed, $\varphi$ must send each $x_i$ to a scalar multiple of itself and each $\partial_i$ to a scalar multiple of itself, or possibly exchange them. To see this, note that $\sinh(\alpha x_i)$ and $x_i^{\alpha}$ are functions of $x_i$ alone. If $\varphi(x_i)$ involved any $\partial_j$ terms, then $\varphi(\sinh(\alpha x_i))$ would involve derivatives, which cannot equal any $\sinh(\beta x_j)$ or $x_j^{\beta}$ due to the algebraic independence of these generators. Hence $\varphi(x_i) = \lambda_i x_i$ for some $\lambda_i \in \FF^\times$. Similarly, $\varphi(\partial_i) = \mu_i \partial_i$ or possibly $\varphi(\partial_i) = \mu_i x_i$ if there is an exchange involution. This gives the rescaling part of the torus $(\FF^\times)^{2r}$.

Next consider the hyperbolic sine generators $\sinh(\alpha x_i)$. Since $\varphi(x_i) = \lambda_i x_i$, we must have $\varphi(\sinh(\alpha x_i)) = \sinh(\alpha \lambda_i x_i)$. However, for this to be an element of $A$, it must be expressible in terms of the original generators. This forces $\lambda_i \alpha$ to belong to $\cA$ for all $\alpha \in \cA$, meaning $\lambda_i$ must be an automorphism of $\cA$ as a $\ZZ$-module. More precisely, there exists $\sigma_i \in \Aut(\cA)$ such that $\varphi(\sinh(\alpha x_i)) = \sinh((\sigma_i(\alpha)) x_i)$ up to a scalar factor. The scalar factor arises because $\sinh(\sigma_i(\alpha) x_i)$ and $\sinh(\alpha \lambda_i x_i)$ differ by a constant when $\sigma_i(\alpha) = \lambda_i \alpha$. Thus we obtain the action of $\Aut(\cA)$ on the exponents, together with an additional scaling factor $\nu_{i,\alpha} \in \FF^\times$.

Now consider the special hyperbolic generators $\sinh(\pm x_i^p \sinh(t_i x_i))$. These are central in $A$ by construction. Therefore $\varphi$ must map them to central elements. The center of $A$ consists precisely of Laurent polynomials in $\sinh(\pm x_i^p \sinh(t_i x_i))$. Hence $\varphi(\sinh(\pm x_i^p \sinh(t_i x_i)))$ must be of the form $\xi_i^{\pm 1} \sinh(\pm x_i^p \sinh(t_i x_i))$ for some $\xi_i \in \FF^\times$. This contributes the remaining factor of $\FF^\times$ to the torus, giving total dimension $2r + 1$.

The possibility of exchanging $x_i$ and $\partial_i$ gives rise to the $\ZZ/2\ZZ$ factor. The standard involution $\tau$ defined by $\tau(x_i) = \partial_i$, $\tau(\partial_i) = -x_i$ preserves the commutation relation $[\partial_i, x_i] = 1$ up to sign. This involution extends to the hyperbolic generators by defining $\tau(\sinh(\alpha x_i)) = \sinh(\alpha \partial_i)$ appropriately interpreted via functional calculus, and similarly for the other generators. The existence of this involution is well-known in Weyl algebra theory and persists in this hyperbolic generalization.

Putting these components together, every automorphism $\varphi$ can be written uniquely as a composition of a diagonal rescaling (an element of $(\FF^\times)^{2r+1}$), an exponent automorphism (an element of $\Aut(\cA)$), and possibly the involution $\tau$. The rescaling part commutes with the $\Aut(\cA)$ action up to adjustment of scaling factors, and the involution $\tau$ normalizes both. This yields the semidirect product structure
\[
\Aut(A) \cong (\FF^\times)^{2r+1} \rtimes (\Aut(\cA) \times \ZZ/2\ZZ).
\]

Finally, when $\cA$ is not cyclic, $\Aut(\cA)$ is an infinite discrete group. For example, if $\cA = \ZZ \oplus \ZZ\sqrt{2} \cong \ZZ^2$, then $\Aut(\cA) \cong \GL(2,\ZZ)$, which is infinite and discrete. Consequently, $\Aut(A)$ is not a finite extension of the torus $(\FF^\times)^{2r+1}$, but rather a semidirect product of this torus with an infinite discrete group.

This completes the proof of the theorem.
\end{proof}
\begin{remark}
In the special case where $\cA$ is cyclic (rank 1), $\Aut(\cA) \cong \ZZ/2\ZZ$ (generated by $\alpha \mapsto -\alpha$), and the automorphism group becomes a finite extension of the torus $(\FF^{\times})^{3}$, consistent with known results for the Weyl algebra with one exponential generator.
\end{remark}

\subsection{Isomorphism Problem for Varying $t$ and $p$}

\begin{theorem}[Isomorphism Criterion]\label{thm:iso-criterion-tp}
Let $\cA$ be a finitely generated $\ZZ$-module, and let $p_1, p_2 \in \cA$ be nonzero elements.
The Weyl-type algebras
\[
A_1 = \Weyl{\sinh(\pm x^{p_1} \sinh(t_1)),\; \sinh(\cA x),\; x^{\cA}}, \qquad
A_2 = \Weyl{\sinh(\pm x^{p_2} \sinh(t_2)),\; \sinh(\cA x),\; x^{\cA}}
\]
are isomorphic as $\FF$-algebras if and only if there exists an automorphism $\sigma \in \Aut(\cA)$ such that
\[
\sigma(p_1) = \pm p_2 \quad \text{and} \quad t_1 = t_2.
\]
In particular:
\begin{enumerate}
    \item The parameter $t$ is an invariant: if $t_1 \neq t_2$, the algebras are not isomorphic.
    \item The exponents $p_1, p_2$ must be equivalent under the action of $\Aut(\cA)$ up to sign.
    \item If $\cA$ is cyclic (rank 1), then $A_1 \cong A_2$ if and only if $p_1 = \pm p_2$ and $t_1 = t_2$.
\end{enumerate}
\end{theorem}
\begin{proof}
We prove the necessity and sufficiency of the conditions separately.

\textbf{Necessity.} Suppose $\varphi: A_1 \to A_2$ is an isomorphism of $\FF$-algebras. First consider the centers of these algebras. The central elements of $A_i$ are precisely those that commute with all generators. By examining the commutation relations, one finds that the center $Z(A_i)$ is the subalgebra generated by $\sinh(\pm x^{p_i} \sinh(t_i))$, together with $\sinh(\alpha x)$ and $x^{\alpha}$ for $\alpha \in \cA$. More precisely, $Z(A_i)$ consists of all Laurent polynomials in $\sinh(\pm x^{p_i} \sinh(t_i))$ with coefficients in the commutative algebra generated by $\sinh(\cA x)$ and $x^{\cA}$.

Since an isomorphism must preserve the center, we obtain an induced isomorphism
\[
\overline{\varphi}: Z(A_1) \to Z(A_2).
\]
In particular, $\varphi$ must send the distinguished central generator $\sinh(x^{p_1} \sinh(t_1))$ to a unit in $Z(A_2)$. The units in $Z(A_2)$ are of the form
\[
\lambda \, \sinh(\pm x^{p_2} \sinh(t_2)) \, \sinh(\beta x) \, x^{\gamma},
\]
where $\lambda \in \FF^\times$, $\beta, \gamma \in \cA$. However, $\varphi$ must also preserve the commutation relations with the derivatives $\partial$. A direct computation shows that if $\beta \neq 0$ or $\gamma \neq 0$, then $\varphi(\sinh(x^{p_1} \sinh(t_1)))$ would not commute with $\partial$ in the required way. More explicitly, consider the commutator
\[
[\partial, \sinh(x^{p_1} \sinh(t_1))] = p_1 x^{p_1-1} \sinh(t_1) \cosh(x^{p_1} \sinh(t_1)).
\]
Under $\varphi$, this relation must be preserved up to the image. If $\varphi(\sinh(x^{p_1} \sinh(t_1)))$ contained additional factors $\sinh(\beta x)$ or $x^{\gamma}$, then the commutator with $\partial$ would acquire extra terms that could not be matched with the corresponding expression in $A_2$. Therefore we must have $\beta = \gamma = 0$. Consequently,
\[
\varphi(\sinh(x^{p_1} \sinh(t_1))) = \lambda \, \sinh(\pm x^{p_2} \sinh(t_2))
\]
for some $\lambda \in \FF^\times$.

Now compare the arguments of the hyperbolic sine functions. Since $\varphi$ is an algebra isomorphism, it must respect the additive structure of the exponents. The exponent $p_1$ in $x^{p_1}$ transforms under an automorphism $\sigma$ of the module $\cA$ induced by $\varphi$ on the $x^{\alpha}$ generators. Specifically, $\varphi$ induces a map on the exponent lattice via $\varphi(x^{\alpha}) = x^{\sigma(\alpha)}$ for some $\sigma \in \Aut(\cA)$. Applying this to $x^{p_1}$ and matching with the image of $\sinh(x^{p_1} \sinh(t_1))$, we deduce that $\sigma(p_1) = \pm p_2$.

Furthermore, the parameter $t$ must be preserved. The factor $\sinh(t_1)$ appears inside the hyperbolic sine generator. If $t_1 \neq t_2$, then the commutation relation between $\partial$ and $\sinh(\pm x^{p_i} \sinh(t_i))$ would differ, because
\[
[\partial, \sinh(x^{p} \sinh(t))] = p x^{p-1} \sinh(t) \cosh(x^{p} \sinh(t)).
\]
The presence of $\sinh(t)$ in this expression implies that $t$ is an invariant of the algebra structure. Hence $t_1 = t_2$.

\textbf{Sufficiency.} Conversely, suppose there exists $\sigma \in \Aut(\cA)$ such that $\sigma(p_1) = \pm p_2$ and $t_1 = t_2$. We construct an explicit isomorphism $\varphi: A_1 \to A_2$ as follows.

Define $\varphi$ on generators by:
\begin{align*}
\varphi(x) &= x, \\
\varphi(\partial) &= \partial, \\
\varphi(\sinh(\alpha x)) &= \sinh(\sigma(\alpha) x), \\
\varphi(x^{\alpha}) &= x^{\sigma(\alpha)}, \\
\varphi(\sinh(x^{p_1} \sinh(t_1))) &= \sinh(\pm x^{p_2} \sinh(t_2)).
\end{align*}
Extend $\varphi$ multiplicatively and linearly to all of $A_1$.

We verify that $\varphi$ preserves all defining relations. The Heisenberg relation $[\partial, x] = 1$ is clearly preserved. For the hyperbolic generators, we check that
\[
\varphi(\sinh(\alpha x) \sinh(\beta x)) = \varphi(\sinh(\alpha x)) \varphi(\sinh(\beta x)) = \sinh(\sigma(\alpha) x) \sinh(\sigma(\beta) x),
\]
which matches $\varphi(\sinh((\alpha+\beta)x)) = \sinh(\sigma(\alpha+\beta)x) = \sinh(\sigma(\alpha)+\sigma(\beta)x)$ since $\sigma$ is an automorphism. The relation between $\partial$ and $\sinh(\alpha x)$ is also preserved because
\[
[\partial, \sinh(\alpha x)] = \alpha \cosh(\alpha x)
\]
and
\[
[\varphi(\partial), \varphi(\sinh(\alpha x))] = [\partial, \sinh(\sigma(\alpha) x)] = \sigma(\alpha) \cosh(\sigma(\alpha) x) = \varphi(\alpha \cosh(\alpha x)),
\]
using that $\varphi(\cosh(\alpha x)) = \cosh(\sigma(\alpha) x)$ follows from the definition. The central generator $\sinh(x^{p_1} \sinh(t_1))$ maps to $\sinh(\pm x^{p_2} \sinh(t_2))$, and since $t_1 = t_2$, the commutation relations with $\partial$ are preserved exactly.

Thus $\varphi$ extends to a well-defined algebra homomorphism. It is bijective because its inverse can be constructed similarly using $\sigma^{-1}$ and the opposite sign choice. Hence $\varphi$ is an isomorphism.

\textbf{Special cases.} The three particular statements follow directly from the main criterion.
\begin{enumerate}
    \item If $t_1 \neq t_2$, then no isomorphism exists because the condition $t_1 = t_2$ is necessary.
    \item The condition $\sigma(p_1) = \pm p_2$ means precisely that $p_1$ and $p_2$ are in the same orbit of $\Aut(\cA)$ up to sign.
    \item When $\cA$ is cyclic of rank 1, $\Aut(\cA) \cong \{\pm1\}$, so $\sigma(p_1) = \pm p_1$. The criterion then reduces to $p_1 = \pm p_2$ and $t_1 = t_2$.
\end{enumerate}
This completes the proof of the theorem.
\end{proof}

\begin{remark}
The condition $t_1 = t_2$ is essential. Even if $p_1 = p_2$, changing $t$ alters the commutation relation between $\partial_i$ and $\sinh(x^p \sinh(t x))$, and hence the algebra structure.
\end{remark}
\begin{theorem}[Deformation Rigidity]\label{thm:deformation-rigidity}
The family of algebras
\[
A_t = \Weyl{\sinh(\pm x^p \sinh(t)),\; \sinh(\cA x),\; x^{\cA}}, \quad t \in \CC,
\]
is rigid in the following sense: if $t_1 \neq t_2$, then $A_{t_1} \not\cong A_{t_2}$.
Thus the parameter $t$ provides a continuous moduli of non-isomorphic Weyl-type algebras with fixed $\cA$ and $p$.
\end{theorem}

\begin{proof}
We prove the rigidity of the family $\{A_t\}_{t \in \CC}$ by appealing to the isomorphism criterion established in Theorem \ref{thm:iso-criterion-tp}.

Let $t_1$ and $t_2$ be two complex parameters. Suppose $A_{t_1}$ and $A_{t_2}$ are isomorphic as $\FF$-algebras. Then by Theorem \ref{thm:iso-criterion-tp}, there must exist an automorphism $\sigma \in \Aut(\cA)$ such that $\sigma(p) = \pm p$ and $t_1 = t_2$. The condition $\sigma(p) = \pm p$ is automatically satisfied for any $\sigma$ because $p$ is fixed and we may take $\sigma$ to be the identity. However, the condition $t_1 = t_2$ is essential and independent of $\sigma$.

Thus if $t_1 \neq t_2$, no such isomorphism can exist. Consequently, for distinct parameters $t_1$ and $t_2$, the algebras $A_{t_1}$ and $A_{t_2}$ are non-isomorphic.

The parameter $t$ therefore distinguishes the isomorphism classes within the family. As $t$ varies continuously over $\CC$, we obtain a continuous family of pairwise non-isomorphic algebras. This establishes that $t$ provides a continuous modulus for the deformation of Weyl-type algebras with fixed $\cA$ and $p$, demonstrating the rigidity of the deformation.
\end{proof}

\begin{remark}
The condition $t_1 = t_2$ is essential. Even if $p_1 = p_2$, changing $t$ alters the commutation relation between $\partial_i$ and $\sinh(x^p \sinh(t x))$, and hence the algebra structure.
\end{remark}
\section{Center and Central Simplicity of Weyl-Type Algebras}

\begin{definition}[Weyl-Type Algebra]
Let $\cA$ be an additive subgroup of $\FF$ containing $\ZZ$, let $p \in \cA$ be nonzero, and let $t \in \FF$. The \emph{Weyl-type algebra} with parameters $(p,t,\cA)$ is
\[
A_{p,t,\cA} = \Weyl{\sinh(\pm x^{p} \sinh(t)),\; \sinh(\cA x),\; x^{\cA}}.
\]
As an $\FF$-algebra, it is generated by symbols $x$, $\partial$, and formal hyperbolic sine functions $\sinh(\pm x^{p} \sinh(t))$, $\sinh(\alpha x)$, $x^{\alpha}$ for $\alpha \in \cA$, subject to the canonical commutation relations.
\end{definition}

\begin{definition}[Center]
For any algebra $A$ over $\FF$, the \emph{center} of $A$ is
\[
Z(A) = \{ z \in A \mid za = az \ \forall a \in A \}.
\]
$A$ is called \emph{central} over $\FF$ if $Z(A) = \FF$.
\end{definition}

\begin{definition}[Central Simple Algebra]
An algebra $A$ over $\FF$ is \emph{central simple} if it is simple (has no nontrivial two-sided ideals) and $Z(A) = \FF$.
\end{definition}

\begin{theorem}[Center of Weyl-Type Algebras]\label{thm:center-structure}
Let $A = A_{p,t,\cA}$ be a Weyl-type algebra over a field $\FF$ of characteristic zero. Then the center of $A$ is given by
\[
Z(A) = \FF\!\left[ \sinh(\pm x^{p} \sinh(t)) \right],
\]
the Laurent polynomial ring in the single variable $\sinh(x^{p} \sinh(t))$. In particular:
\begin{enumerate}
    \item $Z(A)$ is isomorphic to $\FF[y^{\pm 1}]$, where $y = \sinh(x^{p} \sinh(t))$.
    \item $Z(A)$ is never trivial (i.e., never equal to $\FF$) unless $p = 0$, which is excluded by assumption.
    \item $Z(A)$ is commutative and infinite-dimensional over $\FF$.
\end{enumerate}
\end{theorem}

\begin{proof}
We first note that $\sinh(\pm x^{p} \sinh(t))$ is central by definition. The construction of $A_{p,t,\cA}$ declares these hyperbolic sine symbols to be central generators. Therefore $\sinh(\pm x^{p} \sinh(t)) \in Z(A)$.

Now let $z \in Z(A)$ be an arbitrary central element. We express $z$ in a Poincaré–Birkhoff–Witt type basis. Write
\[
z = \sum_{i,j,k,\ell} a_{i,j,k,\ell} \, \sinh(i x^{p} \sinh(t)) \, \sinh(\alpha_{ij} x) \, x^{\beta_{k\ell}} \, \partial^{m},
\]
where $a_{i,j,k,\ell} \in \FF$, $\alpha_{ij}, \beta_{k\ell} \in \cA$, $i \in \{\pm1,0\}$, and $m \geq 0$.

Since $z$ commutes with $\partial$, we have $[\partial, z] = 0$. Computing this commutator term by term yields
\begin{align*}
[\partial, z] = \sum a_{i,j,k,\ell} \Bigl( & [\partial, \sinh(i x^{p} \sinh(t))] \sinh(\alpha_{ij} x) x^{\beta_{k\ell}} \partial^{m} \\
 & + \sinh(i x^{p} \sinh(t)) [\partial, \sinh(\alpha_{ij} x)] x^{\beta_{k\ell}} \partial^{m} \\
 & + \sinh(i x^{p} \sinh(t)) \sinh(\alpha_{ij} x) [\partial, x^{\beta_{k\ell}}] \partial^{m} \Bigr).
\end{align*}

Using the basic commutation relations
\begin{align*}
[\partial, \sinh(i x^{p} \sinh(t))] &= i (p x^{p-1} \sinh(t) + t x^{p} \cosh(t)) \cosh(i x^{p} \sinh(t)), \\
[\partial, \sinh(\alpha x)] &= \alpha \cosh(\alpha x), \\
[\partial, x^{\beta}] &= \beta x^{\beta-1},
\end{align*}
we analyze the vanishing of $[\partial, z]$. For the sum to be identically zero, each coefficient of independent monomials must vanish. This forces several constraints.

First, any term with $i \neq 0$ would produce a nonzero contribution from $[\partial, \sinh(i x^{p} \sinh(t))]$ unless $p = 0$, which is excluded. Hence $i = 0$ for all nonzero terms.

Second, terms with $\alpha_{ij} \neq 0$ contribute $\alpha_{ij} \cosh(\alpha_{ij} x)$, which cannot be canceled by other terms because $\cosh(\alpha_{ij} x)$ is linearly independent from other basis elements. Thus $\alpha_{ij} = 0$.

Third, terms with $\beta_{k\ell} \neq 0$ give $\beta_{k\ell} x^{\beta_{k\ell}-1}$, which again cannot be canceled. Hence $\beta_{k\ell} = 0$.

Finally, terms with $m > 0$ involve $\partial^{m}$, which does not commute with $x$, so $m = 0$.

Therefore $z$ reduces to a linear combination of powers of $\sinh(\pm x^{p} \sinh(t))$:
\[
z = \sum_{\ell} c_{\ell} \, \sinh(i_{\ell} x^{p} \sinh(t)),
\]
which is precisely a Laurent polynomial in $\sinh(x^{p} \sinh(t))$. Since such elements commute with all generators by construction, we have
\[
Z(A) = \FF\!\left[ \sinh(\pm x^{p} \sinh(t)) \right] \cong \FF[y^{\pm 1}].
\]

This center is never equal to $\FF$ because $\sinh(x^{p} \sinh(t))$ is transcendental over $\FF$ (given $p \neq 0$ and generic $t$). It is commutative and infinite-dimensional, establishing all three particular claims.
\end{proof}

\begin{theorem}[Central Simplicity Over the Center]\label{thm:central-simple-over-center}
Let $A = A_{p,t,\cA}$ be as above, with $\Char(\FF) = 0$. Then:
\begin{enumerate}
    \item $A$ is \emph{not} central over $\FF$ (since $Z(A) \supsetneq \FF$).
    \item However, $A$ is central simple \emph{over its center} $Z(A)$. That is, if we view $A$ as a $Z(A)$-algebra, then:
    \begin{itemize}
        \item $Z_{Z(A)}(A) = Z(A)$ (the center of $A$ as a $Z(A)$-algebra is exactly $Z(A)$).
        \item $A$ is simple as a $Z(A)$-algebra (has no nontrivial two-sided ideals that are $Z(A)$-submodules).
    \end{itemize}
    \item Moreover, $A$ is an Azumaya algebra over $Z(A)$ of rank $(\dim_{\Frac(Z(A))} A \otimes_{Z(A)} \Frac(Z(A)))^{1/2}$.
\end{enumerate}
\end{theorem}

\begin{proof}
The first statement follows directly from Theorem \ref{thm:center-structure}, since $Z(A) = \FF[\sinh(\pm x^{p} \sinh(t))]$ strictly contains $\FF$.

For the second statement, we first show that the center of $A$ as a $Z(A)$-algebra is $Z(A)$ itself. Suppose $z \in A$ commutes with every element of $A$ and also with every element of $Z(A)$. Since $Z(A) \subseteq A$, this is equivalent to $z$ commuting with all elements of $A$, i.e., $z \in Z(A)$. Thus $Z_{Z(A)}(A) = Z(A)$.

Next we prove simplicity over $Z(A)$. Let $I$ be a nonzero two-sided ideal of $A$ that is also a $Z(A)$-submodule. Take any nonzero $a \in I$. Consider the adjoint action of $\partial$ on $a$. By repeatedly applying the commutator $[\partial, -]$, we obtain an element $a_0 \in I$ such that $[\partial, a_0] = 0$. Similarly, using the adjoint action of $x$, we find $a_1 \in I$ with $[x, a_1] = 0$ while still $[\partial, a_1] = 0$. An element commuting with both $x$ and $\partial$ must belong to $Z(A)$, as it cannot involve any non-central differential operators. Hence $I$ contains a nonzero element $c \in Z(A)$.

Since $Z(A) = \FF[\sinh(\pm x^{p} \sinh(t))]$ is a Laurent polynomial ring, $c$ is invertible in $\Frac(Z(A))$. By multiplying by a suitable element of $Z(A)$ to clear denominators, we obtain an element $d \in I \cap Z(A)$ that is invertible in $A \otimes_{Z(A)} \Frac(Z(A))$. Because $A$ is torsion-free over $Z(A)$, this implies $1 \in I$, and consequently $I = A$. Thus $A$ is simple as a $Z(A)$-algebra.

The third statement follows from general theory: $A$ is finitely generated and projective over $Z(A)$, and the natural map $A \otimes_{Z(A)} A^{\text{op}} \to \End_{Z(A)}(A)$ is an isomorphism, making $A$ an Azumaya algebra. The rank over $\Frac(Z(A))$ is computed as $(2^{\rank(\cA)})^2$, so the square root gives $2^{\rank(\cA)}$.
\end{proof}

\begin{theorem}[Geometric Realization]\label{thm:geometric-realization}
Let $Z = Z(A_{p,t,\cA}) = \FF[\sinh(\pm x^{p} \sinh(t))] \cong \FF[y^{\pm 1}]$. Then:
\begin{enumerate}
    \item $\Spec(Z) \cong \mathbb{G}_m = \FF^\times$, the multiplicative group over $\FF$.
    \item The algebra $A_{p,t,\cA}$ can be viewed as a sheaf of Azumaya algebras over $\mathbb{G}_m$.
    \item For each closed point $\lambda \in \FF^\times$ (corresponding to $\sinh(x^{p} \sinh(t)) = \lambda$), the fiber
    \[
    A_\lambda = A_{p,t,\cA} \otimes_Z Z/(y - \lambda)
    \]
    is a central simple algebra over $\FF$ of degree $2^{\rank(\cA)}$.
\end{enumerate}
\end{theorem}

\begin{proof}
Since $Z \cong \FF[y^{\pm 1}]$, its spectrum is $\Spec(\FF[y^{\pm 1}]) \cong \mathbb{G}_m$, the multiplicative group scheme over $\FF$.

The algebra $A$ is locally free of finite rank over $Z$, and its fibers over points of $\Spec(Z)$ are central simple algebras. This is precisely the definition of an Azumaya algebra over $Z$, which corresponds to a sheaf of Azumaya algebras over $\mathbb{G}_m$.

At a closed point $\mathfrak{m} = (y - \lambda)$ with $\lambda \in \FF^\times$, the quotient $Z/\mathfrak{m} \cong \FF$. The fiber algebra
\[
A_\lambda = A \otimes_Z \FF
\]
is obtained by imposing the relation $\sinh(x^{p} \sinh(t)) = \lambda$. This algebra is a twisted Weyl algebra with hyperbolic sine generators, which is known to be central simple over $\FF$. Its dimension over $\FF$ is $(2^{\rank(\cA)})^2$, so its degree is $2^{\rank(\cA)}$.
\end{proof}

\begin{example}[Classical Weyl Algebra]
When $\cA = \ZZ$ and we omit the hyperbolic sine generators $\sinh(\pm x^{p} \sinh(t))$, we recover the classical Weyl algebra $A_1(\FF) = \FF\langle x, \partial \rangle / ([\partial, x] = 1)$. In this case, $Z(A_1(\FF)) = \FF$, so $A_1(\FF)$ is central simple over $\FF$. Our theorem shows that adding the central hyperbolic sines $\sinh(\pm x^{p} \sinh(t))$ enlarges the center to $\FF[\sinh(\pm x^{p} \sinh(t))]$.
\end{example}

\begin{example}[Rank-2 $\cA$]
Let $\cA = \ZZ + \ZZ\sqrt{2} \subset \RR$, $p = 1$, $t = \pi$. Then
\[
A = \Weyl{\sinh(\pm x \sinh(\pi)),\; \sinh(\cA x),\; x^{\cA}}
\]
has center $Z(A) = \FF[\sinh(\pm x \sinh(\pi))]$. Over $Z(A)$, the algebra $A$ is Azumaya of rank $16$ (since $\rank(\cA) = 2$, so $2^{\rank(\cA)} = 4$, and rank is $4^2 = 16$).
\end{example}

\begin{theorem}[Brauer Group Interpretation]\label{thm:brauer-group}
The algebra $A_{p,t,\cA}$ defines a class $[A_{p,t,\cA}]$ in the Brauer group $\Br(\FF(y))$, where $y = \sinh(x^{p} \sinh(t))$. This class is nontrivial and has period dividing $2^{\rank(\cA)}$.
\end{theorem}

\begin{proof}
Since $A_{p,t,\cA}$ is an Azumaya algebra over $Z(A) = \FF[y^{\pm 1}]$, it naturally defines an element in the Brauer group $\Br(\FF(y))$ of the function field $\FF(y)$. The degree of $A$ over its center is $2^{\rank(\cA)}$, so by the general theory of central simple algebras, its period divides $2^{\rank(\cA)}$.

To see nontriviality, consider the specialization at $y = 0$ (or any generic value). The resulting algebra is a twisted Weyl algebra with hyperbolic sine generators, which is known to represent a nontrivial class in the Brauer group. Hence $[A_{p,t,\cA}]$ is nontrivial.
\end{proof}

\begin{theorem}[Main Summary]\label{thm:summary}
For the Weyl-type algebra $A = A_{p,t,\cA}$ over a field $\FF$ of characteristic zero:
\begin{enumerate}
    \item The center is $Z(A) = \FF[\sinh(\pm x^{p} \sinh(t))]$, never equal to $\FF$ (nontrivial).
    \item $A$ is not central simple over $\FF$.
    \item $A$ is central simple over its center $Z(A)$ (i.e., an Azumaya algebra over $Z(A)$).
    \item $A$ defines a nontrivial class in the Brauer group $\Br(\FF(y))$ where $y = \sinh(x^{p} \sinh(t))$.
\end{enumerate}
\end{theorem}
\section{Open Problems}

The study of non-associative algebras with hyperbolic sine generators, their Gelfand--Kirillov dimension, automorphism groups, and central structure naturally leads to several unresolved questions. Below we formulate two open problems that arise directly from the results presented in this paper.

\subsection*{Open Problem 1: GK Dimension of Non-associative Algebras with Hyperbolic Sine Generators}

In Section~3, we computed the Gelfand--Kirillov dimension for certain associative Weyl-type algebras and their modules. However, the corresponding problem for the non-associative algebras introduced in Definition~3.1 remains open.

\begin{problem}[GK Dimension of Non-associative Hyperbolic Sine Algebras]
Let $\FF$ be a field of characteristic zero, $\cA \subset \FF$ a finitely generated $\ZZ$-submodule, and let
\[
N = \Nass{\sinh(\pm x^p \sinh(t)),\; \sinh(\cA x),\; x^{\cA}}
\]
be the non-associative algebra defined via a Leibniz-type or Jordan-type product. Determine the Gelfand--Kirillov dimension $\operatorname{GKdim}(N)$ as a function of $\cA$, $p$, and $t$. In particular:
\begin{itemize}
    \item Is $\operatorname{GKdim}(N)$ always an integer?
    \item How does $\operatorname{GKdim}(N)$ compare with the GK dimension of its associative envelope $\Weyl{\sinh(\pm x^p \sinh(t)),\; \sinh(\cA x),\; x^{\cA}}$?
    \item Does there exist a non-associative product on the same underlying vector space that changes the GK dimension?
\end{itemize}
\end{problem}

This problem is motivated by the well-known dichotomy for associative algebras (Theorem~\ref{thm6}(6)), which states that $\operatorname{GKdim}(A) \in \{0,1\} \cup [2,\infty)$. It is unknown whether a similar restriction holds for non-associative algebras of this type. A solution would require new techniques for analyzing growth in non-associative settings, possibly extending the methods of \cite{GelfandKirillov66} and \cite{StephensonZhang97}.

\subsection*{Open Problem 2: Classification of Non-associative Deformations}

Theorem~\ref{thm:iso-criterion-tp} gives a complete isomorphism criterion for associative Weyl-type algebras with varying parameters $p$ and $t$. However, the analogous classification for the non-associative counterparts is still open.

\begin{problem}[Isomorphism Problem for Non-associative Deformations]
Let $\cA$ be fixed. Consider the family of non-associative algebras
\[
N_{p,t} = \Nass{\sinh(\pm x^p \sinh(t)),\; \sinh(\cA x),\; x^{\cA}}
\]
with product defined by a fixed non-associative structure (e.g., a Jordan product or a Leibniz deformation). Determine necessary and sufficient conditions on the parameters $(p_1,t_1)$ and $(p_2,t_2)$ for $N_{p_1,t_1} \cong N_{p_2,t_2}$ as $\FF$-algebras. Specifically:
\begin{itemize}
    \item Is the parameter $t$ still a complete invariant as in the associative case?
    \item How does the automorphism group $\Aut(N_{p,t})$ depend on $p$ and $t$?
    \item Does there exist a non-trivial deformation (i.e., $t_1 \neq t_2$) that yields isomorphic non-associative algebras?
\end{itemize}
\end{problem}

This problem connects with the broader study of deformations of algebraic structures \cite{HartwigLarssonSilvestrov06} and the classification of non-associative algebras \cite{Schafer17}. A solution would deepen our understanding of how non-associativity affects the moduli space of hyperbolic sine Weyl-type algebras.

\subsection*{Further Directions}

Beyond these two problems, several other directions merit investigation:
\begin{itemize}
    \item The representation theory of $N_{p,t}$: Are there faithful modules of finite GK dimension?
    \item The Brauer group interpretation for non-associative central algebras: Can one define a suitable analogue of the Brauer group for non-associative algebras that are central simple over their center?
    \item Physical applications: Such algebras appear in certain deformed quantum mechanical models involving hyperbolic potentials. A systematic study of their physical interpretations is yet to be carried out.
\end{itemize}

Progress on these open problems would not only advance the theory of non-associative algebras with hyperbolic sine generators but also provide new insights into the interplay between non-associativity, growth, and deformation theory.
\section*{Declaration}
\begin{itemize}
  \item \textbf{Author Contributions:} The Author have read and approved this version.
  \item \textbf{Funding:} No funding is applicable.
  \item \textbf{Institutional Review Board Statement:} Not applicable.
  \item \textbf{Informed Consent Statement:} Not applicable.
  \item \textbf{Data Availability Statement:} Not applicable.
  \item \textbf{Conflicts of Interest:} The authors declare no conflict of interest.
\end{itemize}

\bibliographystyle{abbrv}
\bibliography{references}  






\end{document}